\documentclass[10pt,a4paper]{article}

\usepackage{amsmath,amssymb, bbm}
\usepackage{color}
\newcommand{\C}{\mathbb{C}}
\newcommand{\R}{\mathbb{R}}
\newcommand{\N}{\mathbb{N}}
\newcommand{\Z}{\mathbb{Z}}

\newcommand{\B}{\mathbb{B}}

\newcommand{\Bo}{\mathbb{B}_0}
\newcommand{\Co}{\mathbb{C}_0}
\newcommand{\Do}{\mathbb{D}_0}
\def\cA{{\mathcal A}}
\def\cC{{\mathcal C}}
\def\cB{{\mathcal B}}

\def\cF{{\mathcal F}}
\def\cP{{\mathcal P}}
\def\cQ{{\mathcal Q}}
\newcommand{\ee}{\varepsilon}
\renewcommand{\aa}{\alpha}
\newcommand{\aA}{\alpha_1}
\newcommand{\aB}{\alpha_2}
\newcommand{\aC}{\alpha_3}

\renewcommand{\div}{{\rm div}\,}

\newcommand{\Sum}{\displaystyle \sum}
\newcommand{\Hs}{\dot{H^s}}
\newcommand{\bd}{\beta \delta}

\def\d{\partial}
\def\ddl{\dot \Delta_l}
\def\ddj{\dot \Delta_j}
\def\ddq{\dot \Delta_q}
\def\ddk{\dot \Delta_k^v}
\def\tilde{\widetilde}
\def\hat{\widehat}
\newcommand{\D}{\Delta}
\newcommand{\DF}{\Delta_F}

\newcommand{\n}{\nabla}
\newcommand{\G}{\Gamma}

\newcommand{\Fe}{F_{ext}}
\newcommand{\Om}{\Omega}
\newcommand{\Ome}{\Omega_\varepsilon}
\newcommand{\tOm}{\tilde{\Omega}_{QG}}

\newcommand{\ve}{v_\ee}
\newcommand{\Ue}{U_\ee}
\newcommand{\Uqg}{U_{\ee,QG}}
\newcommand{\Uosc}{U_{\ee, osc}}
\newcommand{\Uoe}{U_{0,\ee}}
\newcommand{\Uoqg}{U_{0,\ee,QG}}
\newcommand{\Uoosc}{U_{0,\ee, osc}}
\newcommand{\tUqg}{\tilde{U}_{QG}}
\newcommand{\tvqg}{\tilde{v}_{QG}}
\newcommand{\tUoqg}{\tilde{U}_{0, QG}}
\newcommand{\We}{W_\ee}
\newcommand{\Wet}{W_\ee^T}
\newcommand{\de}{\delta_{\ee}}
\newcommand{\Phie}{\Phi_\ee}
\newcommand{\Thee}{\theta_\ee}

\renewcommand{\Re}{R_\ee}
\newcommand{\re}{r_\ee}
\newcommand{\interv}{[\frac12 -\eta \delta,\frac12 +\eta \delta]}
\newcommand{\intervb}{[\frac12 -\frac32 \delta, \frac12 +\delta]}
\newcommand{\td}{\tilde{\delta}_{n,m}}
\newcommand{\cPrr}{\cP_{\re, \Re}}
\newcommand{\Cdn}{C_{\delta, \nu_0}}

\newtheorem{thm}{Theorem}
\newtheorem{lem}{Lemma}

\newtheorem{prop}{Proposition}
\newtheorem{defi}{Definition}

\newtheorem{rem}{Remark}

\title{Enhanced convergence rates and asymptotics for a dispersive Boussinesq-type system with large ill-prepared data}

\author{Fr\'ed\'eric Charve\footnote{Universit\'e Paris-Est Cr\'eteil, Laboratoire d'Analyse et de Math\'ematiques Appliqu\'ees (UMR 8050), 61 Avenue du G\'en\'eral de Gaulle, 94 010 Cr\'eteil Cedex (France). E-mail: frederic.charve@u-pec.fr}}

\date{}
\begin{document}
%\tableofcontents

\maketitle

\begin{abstract} In this article we obtain, for a stratified, rotating, incompressible Navier-Stokes system, generalized asymptotics as the Rossby number $\ee$ goes to zero (without assumptions on the diffusion coefficients). For ill-prepared, less regular initial data with large blowing-up norm in terms of $\ee$, we show global well-posedness and improved convergence rates (as a power of $\ee$) towards the solution of the limit system, called the 3D quasi-geostrophic system. Aiming for significant improvements required us to avoid as much as possible to resort to classical energy estimates involving oscillations. Our approach relies on the use of structures and symmetries of the limit system, and of highly improved Strichartz-type estimates.
\end{abstract}

\section{Introduction}

\subsection{Geophysical fluids}

The Primitive System (also called Primitive Equations, see for example \cite{Chemin2, BMN5}) is a rotating Boussinesq-type system used to describe geophysical fluids located at the surface of the Earth (in a large physical extent) under the assumption that the vertical motion is much smaller than the horizontal one. Two phenomena exert a crucial influence on geophysical fluids: the Coriolis force induced by the rotation of the Earth around its axis and the vertical stratification of the density induced by gravity. The former induces a vertical rigidity in the fluid velocity as described by the Taylor-Proudman theorem, and the latter induces a horizontal rigidity to the fluid density: heavier masses lay under lighter ones.

In order to measure the importance of these two concurrent phenomena, physicists defined two numbers: the Rossby number $Ro$ and the Froude number $Fr$. We refer to the introduction of \cite{FC, FCpochesLp, FCF1} for more details and to \cite{BeBo, Cushman, Sadourny, Pedlosky} for an in-depth presentation.

The smaller are these numbers, the more important become these two phenomena and we will consider the Primitive Equations in the whole space, under the Boussinesq approximation and when both phenomena are of the same scale i.-e. $Ro=\ee$ and $Fr=\ee F$ with $F>0$. In what follows $\ee$ will be called the Rossby number and $F$ the Froude number. The system is then written as follows (we refer to \cite{Chemin2, BMN5} for the model):
\begin{equation}
\begin{cases}
\d_t \Ue +\ve\cdot \n \Ue -L \Ue +\frac{1}{\ee} \cA \Ue=\frac{1}{\ee} (-\n \Phie, 0),\\
\div \ve=0,\\
{\Ue}_{|t=0}=U_{0,\ee}.
\end{cases}
\label{PE}
\tag{$PE_\ee$}
\end{equation}
The unknowns are $\Ue =(\ve, \Thee)=(\ve^1, \ve^2, \ve^3, \Thee)$ (where $\ve$ denotes the velocity of the fluid and $\Thee$ the scalar potential temperature), and $\Phie$ which is called the geopotential and gathers the pressure term and centrifugal force. The diffusion operator $L$ is defined by
$$
L\Ue \overset{\mbox{def}}{=} (\nu \D \ve, \nu' \D \Thee),
$$
where $\nu, \nu'>0$ are the kinematic viscosity and the thermal diffusivity. The matrix $\cA$ is defined by
$$
\cA \overset{\mbox{def}}{=}\left(
\begin{array}{llll}
0 & -1 & 0 & 0\\
1 & 0 & 0 & 0\\
0 & 0 & 0 & F^{-1}\\
0 & 0 & -F^{-1} & 0
\end{array}
\right).
$$
We will also precise later the properties satisfied by the sequence of initial data $\Uoe$ (as $\ee$ goes to zero). Let us now state some remarks about this system (we refer to the introductions of \cite{FC, FC5, FCpochesLp, FCF1} for more precisions):
\begin{itemize}
 \item This system generalizes the well-known rotating fluids system (see \cite{CDGG, CDGG2, CDGGbook}). The penalized terms (which are are divided by the small parameter $\ee$), namely $\cA \Ue$ and the geopotential, will impose a special structure to the limit when $\ee$ goes to zero.
 \item As $\cA$ is skew-symmetric, and thanks to the incompressibility, any energy method (that is based on $L^2$ or $H^s/\Hs $ inner products) will not "see" these penalized terms and will work as for the classical incompressible Navier-Stokes system ($\cA U\cdot U=0$ and $(\nabla \Phie, v_\ee)_{H^s/\Hs}=0$). Therefore the Leray and Fujita-Kato theorems provide global in time weak solutions if $U_{0,\ee}\in L^2$ and local in time unique strong solutions if $U_{0,\ee} \in \dot{H}^\frac{1}{2}$ (global for small initial data).
 \item There are two distinct regimes wether $F\neq 1$ or $F=1$: the first one features very important dispersive properties. In the second case, the operators are simpler but we cannot rely on Strichartz estimates and the methods are completely different (see \cite{Chemin2, FCF1}). \emph{In the present article we focus on the case $F\neq 1$}.
\end{itemize}

\subsection{Strong solutions}

As explained before, thanks to the skew-symmetry of matrix $\cA$, any computation involving $L^2$ or Sobolev inner-products will be the same as for the Navier-Stokes system. So given the regularity of the initial data (even if some norms can blow up in $\ee$), we can adapt the Leray and Fujita-Kato theorems as well as the classical weak-strong uniqueness results: for a fixed $\ee>0$, if $\Uoe\in \dot{H}^\frac{1}{2}(\R^3)$, we denote as $\Ue$ the unique strong solution of System \eqref{PE}, defined on $[0,T]$ for all $0<T<T_\ee^*$. In addition, if the lifespan $T_\ee^*$ is finite then we have (blow up criterion):
\begin{equation}
 \int_0^{T_\ee^*} \|\n \Ue(\tau)\|_{\dot{H}^\frac{1}{2}(\R^3)}^2 d\tau=\infty.
 \label{critereexpl}
\end{equation}
Moreover, if $\Uoe \in \Hs $ then we also can propagate the regularity as done for the Navier-Stokes system.

\subsection{The limit system, the QG/osc decomposition}
\label{sol}
We are interested in the asymptotics, as the small parameter $\ee$ goes to zero. Let us recall that the limit system is a transport-diffusion system coupled with a Biot-Savart inversion law and is called the 3D quasi-geostrophic system:
\begin{equation}
\begin{cases}
\d_t \tOm +\tvqg .\n \tOm -\G \tOm =0,\\
\tUqg =(\tvqg ,\tilde{\theta}_{QG})=(-\partial_2, \partial_1, 0, -F\partial_3) \DF^{-1} \tOm,
\end{cases}
\label{QG}\tag{$QG$}
\end{equation}
where the operator $\G$ is defined by:
$$
\G \overset{def}{=} \D \DF^{-1} (\nu \d_1^2 +\nu \d_2^2+ \nu' F^2 \d_3^2),
$$
with $\DF=\d_1^2 +\d_2^2 +F^2 \d_3^2$. Moreover we also have the relation
$$
\tOm=\d_1 \tUqg ^2 -\d_2 \tUqg ^1 -F \d_3 \tUqg ^4 =\d_1 \tvqg ^2 -\d_2 \tvqg ^1 -F \d_3 \tilde{\theta}_{QG}.
$$
\begin{rem}
\sl{The operator $\D_F$ is a simple anisotropic Laplacian but $\G$ is in general a tricky non-local diffusion operator of order 2 (except in the case $F=1$ where $\D_F=\D$ and $\G=\nu \d_1^2 +\nu \d_2^2+ \nu' \d_3^2$, or in the case $\nu=\nu'$ where $\G=\nu \Delta$). We refer to \cite{FCestimLp, FCpochesLp} for an in-depth study of $\G$ in the general case (neither its Fourier kernel nor singular integral kernel have a constant sign and no classical result can be used).}
\end{rem}
This limit system is first formally derived then rigourously justified (see \cite{Chemin2, FC}). Led by the limit system we introduce the following decomposition: for any 4-dimensional vector field $U=(v, \theta)$ we define its potential vorticity $\Om(U)$:
$$
\Om(U)\overset{def}{=} \d_1 v^2 -\d_2 v^1 -F\d_3 \theta,
$$
then its quasi-geostrophic and oscillating (or oscillatory) parts:
\begin{equation}
U_{QG}=\cQ (U) \overset{def}{=} \left(
\begin{array}{c}
-\d_2\\
\d_1\\
0\\
-F\d_3
\end{array}
\right) \D_F^{-1} \Om (U), \quad \mbox{and} \quad U_{osc}=\cP (U) \overset{def}{=} U-U_{QG}.
\end{equation}
As emphasized in \cite{FC,FC5} this is an orthogonal decomposition of 4-dimensional vector fields (similar to the Leray orthogonal decomposition into divergence-free and gradient vector fields) and if $\cQ$ and $\cP$ are the associated orthogonal projectors on the quasi-geostrophic or oscillating fields, they satisfy (see \cite{Chemin2, FC, FC2, FCpochesLp, FCF1}):
\begin{prop}
\sl{For any function $U=(v, \theta) \in \Hs $ (for some $s$) we have:
\begin{enumerate}
\item $\cP$ and $\cQ$ are pseudo-differential operators of order 0.
\item For any $s\in\R$, $(\cP(U)|\cQ(U))_{H^s/\Hs} =(\cA U|\cP(U))_{H^s/\Hs}=0$ (when defined).
\item $\cP(U)=U \Longleftrightarrow \cQ(U)=0\Longleftrightarrow \Om(U)=0$.
\item $\cQ(U)=U \Longleftrightarrow \cP(U)=0\Longleftrightarrow$ there exists a scalar function $\Phi$ such that $U=(-\d_2,\d_1,0,-F\d_3) \Phi$. Such a vector field is said to be quasi-geostrophic (or QG) and is also divergence-free.
\item If $U=(v, \theta)$ is a quasi-geostrophic vector field, then $v\cdot \n \Om(U)=\Om(v\cdot \n U)$ and $\G U=\cQ (L U)$.
\item Denoting by $\mathbb{P}$ the Leray orthogonal projector on divergence-free vectorfields, $\mathbb{P} \cP=\cP\mathbb{P}$ and $\mathbb{P} \cQ=\cQ\mathbb{P}=\cQ$.
\end{enumerate}
\label{propdecomposcqg}
}
\end{prop}
Thanks to this, System \eqref{QG} can for example be rewritten into one of the equivalent following velocity formulations:
\begin{equation}
\begin{cases}
\d_t \tUqg +\cQ(\tvqg .\n \tUqg) -\G \tUqg =0,\\
\tUqg =\mathcal{Q} (\tUqg ), \mbox{ (or equivalently } \mathcal{P} (\tUqg )=0),\\
{\tilde{U}_{QG|t=0}= \tUoqg,}
\end{cases}
\label{QG1}\tag{$QG$}
\end{equation}

or

\begin{equation}
\begin{cases}
\d_t \tUqg  +\tvqg .\n \tUqg  -L\tUqg = \cP \tilde{\Phi}_{QG},\\
\tUqg =\mathcal{Q} (\tUqg ),\\
{\tilde{U}_{QG|t=0}= \tUoqg,}
\label{QG2}\tag{$QG$}
\end{cases}
\end{equation}

\begin{rem}\sl{
We recall that Theorem $2$ from \cite{FC2} claims that if $\tilde{U}_{0,QG} \in H^1$ then System \eqref{QG} has a unique global solution $\tUqg  \in \dot{E}^0 \cap \dot{E}^1$ (see below for the space notation). We refer to \cite{FC2, FC4} and to the next sections for more precisions.}
\end{rem}
\begin{rem}
\sl{It is natural to investigate the link between the quasi-geostrophic/oscillating parts decomposition of the initial data and the asymptotics when $\ee$ goes to zero. This leads to the notion of well-prepared/ill-prepared initial data depending on the fact that the initial data is already close or not to the quasi-geostrophic structure, i.-e. when the initial oscillating part is small/large (or going to zero/blowing up as $\ee$ goes to zero). In the present article we consider large and ill-prepared initial data with very large oscillating parts depending on $\ee$.}
\label{illprepared}
\end{rem}
Going back to System \eqref{PE}, we introduce $\Ome=\Om(\Ue)$, $\Uqg=\cQ(\Ue)$ and $\Uosc=\cP(\Ue)$. We showed in \cite{FC} that for an initial data in $L^2$ (independant of $\ee$), the oscillating part $\Uosc$ of a weak global Leray solution $\Ue$, goes to zero in $L_{loc}^2(\R_+,L^q(\R^3))$ ($q\in]2,6[$), and the quasi-geostrophic part $\Uqg$ goes to a solution of System \eqref{QG} (with the QG-part of $U_0$ as initial data). This required the study of System \eqref{systdisp}, and its associated matrix in the Fourier space: as explained in details in Proposition \ref{estimvp} when $\nu\neq \nu'$ there are four distinct eigenvalues (it is necessary to perform frequency truncations to obtain their expression). The first one is explicit but discarded as its associated eigenvector is not divergence-free, the second one is real (and mainly linked to the quasigeostrophic part). The last two ones are non-real and mainly linked to the oscillating part.

Let us denote by $\mathbb{P}_i$ ($i\in\{2,3,4\}$) the associated projectors. When $\nu=\nu'$, many simplifications arise (see Remark \ref{remnunu}). Unfortunately none of these simplifications are true anymore in general (when $\nu\neq \nu'$) but we are able to bound their operator norms and prove that the $\mathbb{P}_2$-part of an oscillating divergence-free vectorfield is small (we refer to \cite{FC, FC3}, see also Proposition \ref{estimvp}).

Moreover we are able to obtain Strichartz estimates for the last two projections $\mathbb{P}_{3+4}$. In \cite{FC} we obtained the following Strichartz estimate upon which depended the main result:
$$
\|\mathbb{P}_{3+4} \mathcal{P}_{r,R}f\|_{L^4L^\infty} \leq C_{r,R} \ee^{\frac14} \left( \|\mathcal{P}_{r,R} f_0\|_{L^2} +\|\mathcal{P}_{r,R} F\|_{L^2}\right).
$$
In \cite{FC2} we focussed on strong solutions. We first proved that if the initial QG-part $U_{0,QG}$ is $H^1$ then the limit system has a unique global solution $\tUqg$. We proved that if $U_{0,osc} \in \dot{H}^\frac12$ then $\Ue$ is global if $\ee$ is small enough. For this we filtered some waves: we constructed a solution $\Wet$ of \eqref{systdisp} with a particular external force term (constructed from $\tUqg$) and proved that $\Ue-\tUqg -\Wet$ goes to zero thanks to a generalization of the previous Strichartz estimates (allowing different regularities for the external force term):
$$
\|\mathbb{P}_{3+4} \mathcal{P}_{r,R}f\|_{L^2 L^\infty} \leq C_{r,R} \ee^{\frac14} \left(\|\mathcal{P}_{r,R} f_0\|_{L^2} +\|\mathcal{P}_{r,R} F^b\|_{L^1 L^2} +\|\mathcal{P}_{r,R} F^l\|_{L^2 L^2}\right).
$$
In \cite{FC3} we generalized the previous result for initial data depending on $\ee$ and with large oscillating part (bounded by $|\ln|\ln \ee||$ in the general case and $|\ln \ee|$ when $\nu=\nu'$) considering frequency truncations $\cPrr$ with radii depending on $\ee$ allowing us to exhibit explicit convergence rates. In this work we distinguished the case $\nu=\nu'$ for which we were able to produce Strichartz estimates \emph{without} frequency truncations in \emph{inhomogeneous} spaces:
$$
\|\We\|_{L^2 B_{\infty,q}^s} \leq C \ee^{\frac18} \left(\|f_0\|_{B_{2,q}^{s+\frac34}}+\|G\|_{L^1 (B_{2,q}^{s+\frac34})}\right)
$$
In the second part of \cite{FC3}, inspired by the work of Dutrifoy about vortex patches in the inviscid case (see \cite{Dutrifoy1}), and by the work of Hmidi for Navier-Stokes vortex patches (see \cite{TH1}), we investigated the case of initial potential vorticity which is a regularized patch, and very large initial oscillating part (regular but bounded by a negative power of $\ee$) when $\nu=\nu'$. This work was recently generalized in the case $\nu\neq\nu'$ in \cite{FCestimLp, FCpochesLp} where we deeply studied the limit quasi-geostrophic operator $\G$ which is non-local and non radial. In this setting, the fact that $\nu\neq \nu'$ highly complicates every computation.

Let us also mention that in \cite{FC4} we obtained global existence when the initial QG-part is only $H^{\frac12+\eta}$. This required real interpolation methods (inspired from \cite{IP}) in order to obtain economic estimates for the limit system (see \eqref{estimQG}). In \cite{FC5} with V.S. Ngo we studied the asymptotics in the case of evanescent viscosities (as a power of $\ee$) and for simplified oscillating initial data (as the initial QG part is zero, the limit is also zero).
\\

Let us now give a survey on other results on this system. In the non-dispersive setting $F=1$ there are few works: let us mention the seminal work of Chemin \cite{Chemin2} (that we recently generalized in \cite{FCF1}) and the work of Iftimie \cite{Dragos4} in the inviscid case.

In \cite{KMY} the authors distinguish the rotation and stratification effects, in the case $\nu=\nu'$ for initial data in $\dot{H}^\frac12 \cap \dot{H}^1$ and for a special condition $\d_2 u_0^1-\d_1 u_0^2=0$ (the initial potential vorticity only depends on the temperature), they obtain existence of a unique global solution to \eqref{PE} in $\mathcal{C}(\R_0, \dot{H}^1)$ for strong enough rotation and stratification. If the initial data is small in $\dot{H}^\frac12$ they manage to obtain that $\n \Ue\in L^2 \dot{H}^\frac12$.

In \cite{LT} Lee and Takada studied global wellposedness in the case of stratification only (no rotationnal effects), when $\nu=\nu'$ and for large initial oscillating part (independant of $\ee$). They first give global existence of a unique mild solution in $L^4(\R_+, \dot{W}^{\frac12,3}(\R^3))$ for large initial oscillating part in $\dot{H}^s$ ($s\in ]\frac12, \frac58]$, there is a kind of smallness condition, see Remark \ref{PetitIMT}) and small QG-part in $\dot{H}^\frac12$. Then they show global well-posedness in the case $s=\frac12$ and for any initial oscillating part and small QG-part, of a unique mild solution in $\mathcal{C}(\R_+, \dot{H}^\frac12) \cap L^4(\R_+, \dot{W}^{\frac12,3}(\R^3))$.

These results are adaptated to the Primitive system in \cite{IMT}. Iwabuchi, Mahalov and Takada focussed on the case $\nu=\nu'$ and obtained (through stationnary phase methods) the following Strichartz estimates that we state with our notations:
\begin{prop} (\cite{IMT} Theorem 1.1 and Corollary 1.2)
\sl{Assume $F\neq 1$. If $r\in]2,4[$ and $p\in]2,\infty[ \cap [\frac1{2(\frac12-\frac1{r})}, \frac2{3(\frac12-\frac1{r})}]$, there exists a constant $C=C_{F,\nu,p,r}$ such that if $f$ solves the homogeneous \eqref{systdisp},
$$
\|f\|_{L^p(\R_+, L^r)} \leq C \ee^{\frac1{p} -\frac32(\frac12 -\frac1{r})} \|f_0\|_{L^2}.
$$
If $s\in]\frac12, \frac58]$, there exists a constant $C=C(F,s,\nu)$ such that:
$$
\|f\|_{L^4(\R_+, \dot{W}^{s,\frac6{1+2s}})} \leq C \ee^{\frac12(s-\frac12)} \|f_0\|_{\dot{H}^s}.
$$
}
\label{StriIMT}
\end{prop}
From this they are able to obtain through a fixed point argument the following global well-posedness results for initial data (independant of $\ee$) with small quasi-geostrophic part (assume $\nu=\nu'$ and $F\neq 1$):
\begin{itemize}
 \item If $s\in]\frac12, \frac58]$, there exist $\delta_1,\delta_2>0$ (depending on $\nu,F,s$) such that for any $\ee>0$ and any initial data $U_0= U_{0,QG} + U_{0,osc}$ with $(U_{0,QG}, U_{0,osc})\in \dot{H}^\frac12 \times \dot{H}^s$ and
\begin{equation}
  \begin{cases}
   \|U_{0,QG}\|_{\dot{H}^\frac12} \leq \delta_1,\\
  \|U_{0,osc}\|_{\dot{H}^s} \leq \delta_2 \ee^{-\frac12(s-\frac12)},
 \end{cases}
 \label{IMTpetitesse}
\end{equation}
 there exists a unique global mild solution in $L^4(\R_+,\dot{W}^{\frac12,3}(\R^3))$.
 \item There exists $\delta>0$ such that for any initial data $U_0= U_{0,QG} + U_{0,osc}\in \dot{H}^\frac12$ with $\|U_{0,QG}\|_{\dot{H}^\frac12} \leq \delta$, there exists $\ee_0>0$ such that for any $0<\ee<\ee_0$, System \eqref{PE} has a unique global mild solution in $\mathcal{C}(\R_+, \dot{H}^\frac12) \cap L^4(\R_+,\dot{W}^{\frac12, 3}(\R^3))$.
\end{itemize}

Let us also mention works in the periodic case where resonences have to be studied (see for example \cite{IG1, VSN, VSNS1, Scro}), in the rotating fluids system case (see \cite{CDGG, CDGG2, CDGGbook, GIMS, HS, KLT}) or in the inviscid case (see \cite{Dutrifoy1, Dutrifoy2, KLT2, T1, Wild}).
\\

In the present article we wish to generalize our results from \cite{FC2, FC3, FC4} and motivated by the very interesting results in \cite{IMT} we want to obtain full asymptotics (as in \cite{FC5, FCpochesLp}) for very large ill-prepared initial data (less regular, depending on $\ee$ and bounded by a negative power of $\ee$). In our work we will provide global well-posedness results but also precise convergence rates as $\ee$ goes to zero. We also generalize \cite{IMT} in the sense that we consider initial data with \emph{large quasi-geostrophic part} (with low frequencies assumptions) and provide solutions in \emph{homogeneous} energy spaces $\dot{E}^s$ both in the particular case $\nu=\nu'$ \emph{and} in the general case $\nu\neq \nu'$. Let us also mention that our methods closely rely on the special structures and properties of the 3D quasi-geostrophic system.

\subsubsection{Statement of the results}
We will consider general ill-prepared initial data $\Uoe =\Uoosc +\Uoqg$, whose $QG$-part converges to some $\tUoqg$ (without any smallness condition), and whose oscillating part is very large (see below for precisions).
\\

The aim of the present article is to generalize Theorem 3 from \cite{FC2}, Theorems 1.2 and 1.3 from \cite{FC3} and Theorem 4 from \cite{FC4} with the least possible extra regularity for the initial data and the biggest possible blowing-up initial oscillatory part (as a negative power of $\ee$). The energy methods used in \cite{FC, FC2, FC3} would only allow at best an initial blow-up of $\Uoosc$ as $|\ln \ee|^\beta$. Indeed, these methods require the use of energy estimates for the oscillations $\We$ or $\Wet$ and produce large terms involving $exp(\|\Uoosc\|^2)$ that can only be balanced thanks to $\ee^\gamma$ provided by the Strichartz estimates. We need to change our point of view and try to not resort to energy estimates for these oscillations. This will require us to make more flexible dispersive estimates so that the oscillations can be estimated with minimal use of their energy (the only term where it was unavoidable is $F_8$, see below for details). We will here state only the new results. Let us define (in the whole space $\R^3$) the family of spaces $\dot{E}_T^s$ for $s\in \R$,
$$
\dot{E}_T^s=\mathcal{C}_T(\Hs ) \cap L_T^2(\dot{H}^{s+1}),
$$
endowed with the following norm (where we define $\nu_0=\min(\nu, \nu')$, see the appendix for the other notations):
$$
\|f\|_{\dot{E}_T^s}^2 \overset{def}{=}\|f\|_{L_T^\infty \Hs }^2+\nu_0 \int_0^T \|f(\tau)\|_{\dot{H}^{s+1}}^2 d\tau.
$$
When $T=\infty$ we denote $\dot{E}^s$ and the corresponding norm is over $\R_+$ in time. Let us now state the main result of this article (we \emph{do not} assume $\nu=\nu'$).

\begin{thm}
 \sl{Assume $F\neq 1$. For any $\Co\geq 1$, $\delta\in]0,\frac1{10}]$, $\alpha_0>0$, there exist five constants $\varepsilon_0,\eta,\Bo, \kappa, \beta>0$ (depending on $F, \nu, \nu',\Co, \alpha_0$) such that for all $\ee\in]0,\ee_0]$ and all divergence-free initial data $\Uoe= \Uoqg + \Uoosc$ satisfying:
 \begin{enumerate}
  \item $\Uoqg$ converges towards some quasi-geostrophic vectorfield $\tUoqg\in H^{\frac12 + \delta}$ with:
  \begin{equation}
   \begin{cases}
    \|\Uoqg-\tUoqg\|_{H^{\frac12 + \delta}}\leq \C_0 \ee^{\aa_0},\\
    \|\tUoqg\|_{H^{\frac12 + \delta}}\leq \C_0.
   \end{cases}
  \end{equation}
\item $\|\Uoosc\|_{\dot{F}_\delta} \leq \C_0 \ee^{-\kappa \delta}$ where the space $\dot{F}_\delta$ is defined as follows ($q=\frac{2}{1+\delta}$):
$$
\dot{F}_\delta =
\begin{cases}
 \dot{H}^{\frac12 -\delta} \cap \dot{H}^{\frac12 +\delta} & \mbox{ if }\nu=\nu',\\
 \dot{B}_{q,q}^{\frac12} \cap \dot{H}^{\frac12 + \delta} & \mbox{ if }\nu\neq\nu',
\end{cases}
$$
 \end{enumerate}
then System \eqref{QG3} has a unique global solution $\tUqg\in \dot{E}^0 \cap \dot{E}^{\frac12 +\delta}$, and System \eqref{PE} has a unique global solution $\Ue \in \dot{E^s}$ for all $s\in[\frac12-\eta \delta, \frac12+\eta \delta]$, which converges towards $\tUqg$ with the following estimate:
$$
\|\Ue-\tUqg\|_{L^2 L^\infty} \leq \B_0 \ee^{\min(\aa_0, \delta \beta)}.
$$
\label{Th0}
 }
\end{thm}

\begin{rem}
 \sl{In the general case, $\kappa$ is small (less than $\frac1{4000}$), whereas in the case $\nu=\nu'$, $\kappa<\frac12$ (and as close to $\frac12$ as we want). We refer to the next section for a more precise statement of this theorem.}
\end{rem}
\begin{rem}
 \sl{It is interesting to adapt these results to the case with only stratification.}
\end{rem}

\subsection{Precise statement of the main results}

This section is devoted to give the precise statement of Theorem \ref{Th0}, which will be split into two formulations wether we have $\nu=\nu'$ or $\nu\neq \nu'$. This statement requires us to introduce auxilliary systems, which is the object of the first two subsections, and state additional regularity properties for the solution of the limit system (we refer to the third subsection). Then we will state the results we will prove in this article.

\subsubsection{Auxiliary systems in the general case $\nu\neq \nu'$}

\begin{rem}
 \sl{In what follows, we will systematically write, for $f:\R^3 \rightarrow \R^4$, $f\cdot \nabla f=\sum_{i=1}^3 f_i \d_i f$.}
\end{rem}
Following \cite{FC2} we rewrite the primitive system, projecting onto the divergence-free vectorfields ($\mathbb{P}$ is the classical Leray projector):
\begin{equation}
\begin{cases}
{\d_t \Ue -L \Ue+\frac{1}{\varepsilon} \mathbb{P} \mathcal{A} \Ue =-\mathbb{P}(\Ue.\nabla \Ue).} \\
{{\Ue}_{|t=0}=\Uoe.}
\end{cases}
\label{pe2}
\end{equation}
Notice that we can rewrite \eqref{QG} as follows (we also refer to \cite{FC2} where it was first used):
\begin{equation}
\begin{cases}
{\d_t \tUqg -L \tUqg+\frac{1}{\ee} \mathbb{P} \mathcal{A} \tUqg=-\mathbb{P}(\tUqg.\nabla \tUqg) +G,}\\
{\tilde{U}_{QG|t=0}= \tUoqg.}
\end{cases}
\label{QG3}\tag{$QG$}
\end{equation}
where
\begin{equation}
G= G^b + G^l \overset{def}{=}
\mathbb{P} \mathcal{P} (\tUqg. \nabla \tUqg) -F(\nu-\nu')\Delta \Delta_F^{-2}
\operatorname{} \left(
  \begin{array}{c}
-F \partial_2 \partial_3^2 \\
F \partial_1 \partial_3^2 \\
0\\
(\partial_1^2+\partial_2^2)\partial_3
\end{array} \right) \tOm.
\label{G}
\end{equation}

\begin{rem}
 \sl{It is important to notice that $G$ is the sum of two terms, both divergence-free and whose potential vorticity is zero, which is crucial to fully take advantage of \eqref{estimvp2}. We refer to \cite{FC2, FC4} for more details}.
 \label{Gomeganul}
\end{rem}

As explained in \cite{FC, FC2, FC3, FC4, FC5, FCpochesLp}, in the case $F\neq 1$ the oscillatory part enjoys dispersive properties that allow us to obtain Strichartz-type estimates. More precisely the oscillatory part satisfies System \eqref{systdisp} (we refer to the appendix for details), and 
in all the cited articles, we used that the frequency truncated third and fourth projections of the oscillatory part satisfy Strichartz-type estimates as given by Proposition \ref{Strichartz1}. As in \cite{FC2, FC4, FC5}, in the present article we will consider some particular oscillatory terms whose existence is only devoted to absorb some constant terms in order to get the desired convergence rate for the asymptotics as $\ee$ goes to zero.

More precisely, we introduce the following linear system (we refer to the appendix for the notations $\re, \Re$ and $\cPrr$):
\begin{equation}
 \begin{cases}
  \d_t \Wet -L \Wet +\frac{1}{\varepsilon} \mathbb{P} \mathcal{A} \Wet =-\cPrr \mathbb{P}_{3+4} G,\\
  {\Wet}_{|t=0}=\cPrr \mathbb{P}_{3+4} \Uoosc
 \end{cases}
\label{WT}
 \end{equation}
\begin{rem}
 \sl{We recall that it would be useless to consider the free system: indeed the system satisfied by $\Ue-\tUqg$ features $G$ as an external force term which is independant of $\ee$ and blocks any convergence. It is then necessary to absorb a large part of this term which is the reason why we considered such an external force term in System \eqref{WT}. In other words, $\Wet$ is small due to dispersive properties, but still it allows us to "eat" a large part of $G$. We refer to \cite{FC2} for more details.}
\end{rem}

Finally we define $\de=\Ue -\tUqg -\Wet$, which satisfies the following system (see \cite{FC2} for details):
\begin{equation}
 \begin{cases}
  \d_t \de -L\de +\frac{1}{\ee} \mathbb{P} \mathcal{A} \de = \Sum_{i=1}^8 F_i +f^b +f^l,\\
  {\de}_{|t=0}= (\Uoqg-\tUoqg) +(Id-\cPrr)\Uoosc +\cPrr \mathbb{P}_2 \Uoosc,
 \end{cases}
\label{DE}
\end{equation}
where we define:
\begin{equation}
\begin{cases}
 F_1 \overset{def}{=}-\mathbb{P}(\de \cdot \nabla \de), \quad F_2 \overset{def}{=}-\mathbb{P}(\de \cdot \nabla \tUqg), \quad F_3 \overset{def}{=}-\mathbb{P}(\tUqg \cdot \nabla \de),\\
 F_4 \overset{def}{=}-\mathbb{P}(\de \cdot \nabla \Wet),\quad F_5 \overset{def}{=}-\mathbb{P}(\Wet \cdot \nabla \de),\quad F_6 \overset{def}{=}-\mathbb{P}(\tUqg \cdot \nabla \Wet),\\
 F_7 \overset{def}{=}-\mathbb{P}(\Wet \cdot \nabla \tUqg),\quad F_8 \overset{def}{=}-\mathbb{P}(\Wet \cdot \nabla \Wet),\\
 f^b \overset{def}{=} -(Id-\cPrr)G^b -\cPrr \mathbb{P}_2 G^b,\\
 f^l \overset{def}{=} -(Id-\cPrr)G^l -\cPrr \mathbb{P}_2 G^l.  
\end{cases}
 \label{f1f2}
\end{equation}

\subsubsection{Auxiliary systems in the special case $\nu=\nu'$}

In this case, many simplifications arise in the computations of the eigenvalues and eigenvectors of System \eqref{systdisp} (see Remark \ref{remnunu}). In this case, as used in the first part of \cite{FC3}, we can use the following system instead of \eqref{WT}:

\begin{equation}
 \begin{cases}
  \d_t \We -L \We +\frac{1}{\varepsilon} \mathbb{P} \mathcal{A} \We = -G^b,\\
  {\We}_{|t=0}=\Uoosc
 \end{cases}
\label{We}
 \end{equation}

We will be able in the present article to provide for this system much more accurate Strichartz estimates without any frequency restrictions (generalizing the ones obtained in \cite{FC3}). If we denote $\de=\Ue -\tUqg -\We$, it satisfies the following system:

\begin{equation}
 \begin{cases}
  \d_t \de -L\de +\frac{1}{\ee} \mathbb{P} \mathcal{A} \de = \Sum_{i=1}^8 F_i,\\
  {\de}_{|t=0}= \Uoqg-\tUoqg,
 \end{cases}
\label{GE}
\end{equation}

\begin{rem}
 \sl{We choose here to use the same notations as in the general case, the only difference is that $\Wet$ has to be replaced by $\We$.}
\end{rem}

\subsubsection{The limit system}

Let us recall that Theorem $2$ from \cite{FC2} states that when the initial data $\tUoqg$ is in the \emph{inhomogeneous} Sobolev space $H^1$ then System \eqref{QG3} has a unique global solution $\tUqg \in \dot{E}^0\cap \dot{E}^1$, moreover there exists a constant $C=C(\delta)>0$ such that for all $s\in [0,1]$ and all $t\in \R_+$ (and denoting as usual $\nu_0=\min(\nu, \nu')>0$):
$$
\|\tUqg \|_{L_t^\infty \Hs }^2+\nu_0 \int_0^t \|\n\tUqg (\tau)\|_{\Hs }^2 d\tau \leq C(\|\tilde{U}_{0,QG}\|_{L^2}^{1-s}\|\tilde{U}_{0,QG}\|_{\dot{H}^1}^s)^2 \leq C \|\tilde{U}_{0,QG}\|_{H^1}^2.
$$
In \cite{FC4} we used real interpolation methods from Gallagher and Planchon in \cite{IP} (we also refer to the work of C\`alderon in \cite{Calderon}) to obtain a much more accurate estimate, which allowed to bound the energy in $\dot{E}^0\cap \dot{E}^{\frac12 + \delta}$ only with the $H^{\frac12 + \delta}$ initial norm instead of the full $H^1$ norm (we refer to Lemma $2.1$ in \cite{FC4}, our aim was to consider less regular initial data): for any $\delta>0$ there exists a constant $C=\Cdn>0$ such that for all $t\in \R_+$:
\begin{equation}
 \|\tUqg \|_{L_t^\infty H^{\frac12 + \delta}}^2+\nu_0 \int_0^t \|\n \tUqg (\tau)\|_{H^{\frac12 + \delta}}^2 d\tau \leq \Cdn \|\tUoqg\|_{H^{\frac12 + \delta}}^2 \max(1, \|\tUoqg\|_{H^{\frac12 + \delta}}^{\frac{1}{\delta}}).
 \label{estimQG}
\end{equation}
\begin{rem}
 \sl{The reader may wonder why the right-hand side is not simply $\Cdn \|\tUoqg\|_{H^{\frac12 + \delta}}^{2+\frac1{\delta}}$ as stated in \cite{FC4, IP}. This is the right formulation when $\|\tUoqg\|_{H^{\frac12 + \delta}}$ is large (in \cite{FC4} we implicitely focussed on large initial QG part). When it is small, the right-hand side is even simpler: $\Cdn \|\tUoqg\|_{H^{\frac12 + \delta}}^2$.
In the proof in \cite{FC4} of \eqref{estimQG} it is crucial to use Lemma 4.3 from \cite{IP}, and for this, some threshold $j_0\geq 1$ has to be defined:
\begin{itemize}
 \item Either $\|\tUoqg\|_{H^{\frac12 + \delta}}> \frac23 c\nu_0 2^{2\delta}$ (where $c$ is the smallness constant from the Fujita-Kato theorem), and we can define the threshold $j_0$  as stated in \cite{FC4} so that the right-hand side of \eqref{estimQG} is $C_0 (1-2^{-4\delta})^{-2} \left(\frac{3}{2c\nu_0}\right)^{\frac1{\delta}}\|\tUoqg\|_{H^{\frac12 + \delta}}^{2+\frac1{\delta}}$ ($C_0$ is a universal constant).
 \item Or $\|\tUoqg\|_{H^{\frac12 + \delta}}\leq \frac23 c \nu_0 2^{2\delta}$ and then we can simply choose the threshold $j_0=1$ and obtain \eqref{estimQG} with right-hand side that can be simplified into $C_0 (1-2^{-4\delta})^{-2} \|\tUoqg\|_{H^{\frac12 + \delta}}^2$.
\end{itemize}
In other words, the right-hand side of \eqref{estimQG} is in general:
$$
C_0 \big( 1-2^{-4 \delta}\big)^{-2\delta} \|\tUoqg\|_{H^{\frac12 + \delta}}^2 \max \left(1, \frac14 \left(\frac{3}{2c\nu_0}\|\tUoqg\|_{H^{\frac12 + \delta}}\right)^{\frac{1}{\delta}}\right).
$$
}
\end{rem}
Our first result is devoted to the limit system and generalises Theorem $2$ from \cite{FC2} using the precise estimates obtained in \cite{FC4}:

\begin{thm}
 \sl{Let $\delta>0$ and $\tUoqg\in H^{\frac12 + \delta}$ a quasigeostrophic vectorfield (that is $\tUoqg=\cQ \tUoqg$). Then System \eqref{QG3} has a unique global solution in $E^{\frac12 + \delta}=\dot{E}^0\cap \dot{E}^{\frac12 + \delta}$ and the previous estimates holds true.
\label{Th1}
 }
\end{thm}

\subsubsection{Statement in the case $\nu=\nu'$}

\begin{thm}
 \sl{Assume $F\neq 1$. For any $\C_0\geq 1$, $\delta\in]0,\frac{1}{10}]$, $\gamma\in]0,\frac{\delta}{2}[$ and any $\aa_0>0$, if we define $\eta>0$ such that
 $$
 \gamma=(1-2\eta)\frac{\delta}{2} \quad (\mbox{that is } \eta=\frac12(1-\frac{2\gamma}{\delta})),
 $$
 there exists $\ee_0,\B_0>0$ (all of them depending on $F, \nu,\Co, \delta, \gamma, \alpha_0$) such that for all $\ee\in]0,\ee_0]$ and all divergence-free initial data $\Uoe= \Uoqg + \Uoosc$ satisfying:
 \begin{enumerate}
  \item There exists a quasi-geostrophic vectorfield $\tUoqg\in H^{\frac12 + \delta}$ such that
  \begin{equation}
   \begin{cases}
    \|\Uoqg-\tUoqg\|_{H^{\frac12 + \delta}}\leq \C_0 \ee^{\aa_0},\\
    \|\tUoqg\|_{H^{\frac12 + \delta}}\leq \C_0.
   \end{cases}
  \end{equation}
\item $\Uoosc\in \dot{H}^{\frac12} \cap \dot{H}^{\frac12 + \delta}$ with $\|\Uoosc\|_{\dot{H}^{\frac12} \cap \dot{H}^{\frac12 + \delta}} \leq \C_0 \ee^{-\gamma}$,
 \end{enumerate}
then System \eqref{PE} has a unique global solution $\Ue \in \dot{E^s}$ for all $s\in[\frac12, \frac12+\eta \delta]$, and if we define
\begin{itemize}
 \item $\tUqg$ as the unique global solution of \eqref{QG3} in $\dot{E}^0 \cap \dot{E}^{\frac12 +\delta}$,
\item $\We$ as the unique global solution of \eqref{We} in $\dot{E}^{\frac12} \cap \dot{E}^{\frac12 +\delta}$,
\item $\de=\Ue-\tUqg-\We$,
\end{itemize}
then for all $s\in[\frac12, \frac12+\eta \delta]$
\begin{equation}
 \|\de\|_{\dot{E}^s}\leq \B_0 \ee^{\min(\aa_0, \frac{\delta \eta}{2})}.
 \label{estimTh1}
\end{equation}
Moreover if we ask for more low frequency regularity for the initial oscillating part, that is $\Uoosc\in \dot{H}^{\frac12-\delta} \cap \dot{H}^{\frac12 + \delta}$ with $\|\Uoosc\|_{\dot{H}^{\frac12 -\delta} \cap \dot{H}^{\frac12 + \delta}} \leq \C_0 \ee^{-\gamma}$ then \eqref{estimTh1} is true for all $s\in[\frac12 -\eta \delta, \frac12+\eta \delta]$ and we also can get rid of the oscillations $\We$ and obtain that:
$$
\|\Ue-\tUqg\|_{L^2 L^\infty} \leq \B_0 \ee^{\min(\aa_0, \frac{\delta \eta}{2})}.
$$
\label{Th2}
 }
\end{thm}
\begin{rem}
 \sl{Compared to Theorem $1.3$ from \cite{FC3} we highly reduced the regularity of the initial data, only the quasi-geostrophic part lies in a inhomogeneous space and we allow a far greater blowup in $\ee$ for the oscillating part, keeping a satisfying convergence rate as a power of $\ee$ (in accordance with Physicists) for any size of the initial quasi-geostrophic part.}
\end{rem}
%\begin{rem}
% \sl{The previous results can be formulated in two ways: either fixing the size of data $\C_0$ and then obtaining the results for all $\ee\leq \ee_0(\C_0)$ (as we did here and in \cite{FC} to \cite{FCF1}), or stating the results for any $\ee$ provided that the size of the data $\C_0(\ee)$ is small enough as for example in \cite{Dutrifoy1} or \cite{IMT} (see Remark $1.4$).}
%\end{rem}
\begin{rem}
 \sl{Note that in \cite{LT, IMT} there is a smallness condition for the initial quasi-geostrophic part (and also for the oscillating part in some sense). Their result states there exist $\delta_{1,2}>0$ such that for any initial data satisfying \eqref{IMTpetitesse}, there exists a global unique mild solution for any $\ee>0$. This has to be compared with our formulation, where we prove that for any size $\Co$ and any initial data with $\|\Uoqg\|\leq \Co$ and $\|\Uoosc\|\leq \Co \ee^{-\gamma}$, there exists a unique global solution when $\ee\leq \ee_0$.}
 \label{PetitIMT}
\end{rem}

\begin{rem}
 \sl{Compared to the assumptions in \cite{IMT} (Theorems 1.3 and 1.5), we reach the same regularity for the oscillating part, we ask more regularity to the initial QG-part, and we ask more low frequency regularity for both of them (we have to assume $\Uoe\in \dot{H}^{\frac12}$ as we need to consider Fujita-Kato strong solutions):
 $$
 \begin{cases}
  \Uoosc \in \dot{H}^{\frac12} \cap \dot{H}^{\frac12 + \delta} \quad (\dot{H}^{\frac12 + \delta}\mbox{ in } \cite{IMT}),\\
  \Uoqg \in H^{\frac12 + \delta} \quad (\dot{H}^{\frac12}\mbox{ in } \cite{IMT}),
 \end{cases}
 $$
 but \emph{we do not ask any smallness to the initial quasi-geostrophic part}, and we provide global strong solutions \emph{in the energy spaces} $\dot{E}^s$ for any $s\in \interv$ (compared to mild solutions in $L^4(\R_+, \dot{W}^{\frac12,3})$).}
\end{rem}
\begin{rem}
 \sl{At first sight our blow-up rate seems slightly less general than the one from \cite{IMT} (in \cite{IMT} they ask $\ee^{\frac{\delta}2}\|\Uoosc\|_{\dot{H}^{\frac12+\delta}}$ smaller than some $\delta_2>0$, and in the present work, we choose any $\Co$ and ask $\ee^{\gamma}\|\Uoosc\|_{\dot{H}^\frac12 \cap \dot{H}^{\frac12+\delta}} \leq \Co$ for any $\gamma<\frac{\delta}2$) but in our result we look for explicit rates of convergence as powers of $\ee$. We refer to Remark \ref{RqIMT} for more details.}
\end{rem}
\begin{rem}
 \sl{We refer to Remark \ref{RqStri} for a comparision of the Strichartz estimates we use and the ones from \cite{IMT}.}
\end{rem}

\subsubsection{Statement in the general case $\nu\neq\nu'$}

\begin{thm}
 \sl{Assume $F\neq 1$. Let $\delta\in]0,\frac12]$, $q=\frac{2}{1+\delta}$, $\aa_0>0$, $m\in]0, \frac{1}{100}]$, and $M,\eta>0$ such that
 $$
 0<2\eta\leq \frac{M}{m} \leq \frac12 \frac1{5+\delta},
 $$
 let $\gamma_0 \in ]0,\frac{M \delta}{4}]$. If we define $\Re =\ee^{-M}$ and $\re =\ee^m$ then for all $\Co\geq 1$, there exist $\ee_0$, $\Bo$ (all of them depending on $F, \nu, \nu',\Co, \delta, \gamma, \alpha_0$) such that for all initial data $\Uoe=\Uoosc+ \Uoqg$ satisfying:
 \begin{enumerate}
  \item There exists a quasi-geostrophic vectorfield $\tUoqg\in H^{\frac12 + \delta}$ such that
  \begin{equation}
   \begin{cases}
    \|\Uoqg-\tUoqg\|_{H^{\frac12 + \delta}}\leq \C_0 \ee^{\aa_0},\\
    \|\tUoqg\|_{H^{\frac12 + \delta}}\leq \C_0.
   \end{cases}
  \end{equation}
\item $\Uoosc\in \dot{B}_{q,q}^{\frac12} \cap \dot{H}^{\frac12 + \delta}$ with $\|\Uoosc\|_{\dot{B}_{q,q}^{\frac12} \cap \dot{H}^{\frac12 + \delta}} \leq \C_0 \ee^{-\gamma}$,
 \end{enumerate}
then System \eqref{PE} has a unique global solution $\Ue \in \dot{E}^s$ for all $s\in \interv$. Moreover, with the same notations as in Theorem \ref{Th1} (replacing $\We$ by $\Wet$, which involves $m,M$),
\begin{equation}
 \|\de\|_{\dot{E}^s} \leq \Bo \ee^{\min(\aa_0, \frac{M\delta}{4})},
\label{CVdelta}
\end{equation}
and finally, thanks to the Strichartz estimates, we can get rid of the oscillations $\Wet$ and obtain:
$$
\|\Ue-\tUqg\|_{L^2 (\R_+, L^\infty)} \leq \Bo \ee^{\min(\aa_0, \frac{M\delta}{4})}.
$$
\label{Th3}
 }
\end{thm}
\begin{rem}
 \sl{This generalizes the first result from \cite{FC3}: in the present work we reduced the assumptions on high and low frequencies for the initial oscillating part and the choice for $\re$ and $\Re$ now correctly fits the power of $\ee$ provided by the Strichartz estimates, which produces a convergence rate as a power of $\ee$ without any assumption on the viscosities.}
\end{rem}
\begin{rem}
 \sl{The low-frequencies assumption $\Uoosc \in \dot{B}_{q,q}^{\frac{1}{2}}$ is mainly needed to produce a positive power of $\ee$ when estimating $\|\chi (\frac{|D|}{\Re}) \chi (\frac{|D_3|}{\re}) \Uoosc\|_{\Hs }$ (the other need is to reach regularities less than $\frac12$), and the high-frequencies assumption $\Uoosc \in \dot{H}^{\frac{1}{2}+\delta}$ helps to estimate $\|(1-\chi (\frac{|D|}{\Re})) \Uoosc\|_{\Hs }$.
 }
\end{rem}
\begin{rem}
 \sl{The classical Bernstein estimates ensures that $\dot{B}_{q,q}^{\frac{1}{2}} \hookrightarrow \dot{H}^{\frac12 -\frac32 \delta}$ so that $\Uoosc \in \Hs $ for all $s\in \intervb$.}
\end{rem}

The rest of this article is structured as follows: we will first prove Theorem \ref{Th1}, then turn to the proof of Theorem \ref{Th2} in the case $\nu=\nu'$ (much easier computations to obtain the eigenvalues and vectors, but needs more careful use for the Strichartz estimates as $\We$ is not frequency truncated) and we will finish with the proof of Theorem \ref{Th3} (the eigenvectors are not mutually orthogonal anymore, and care is needed for the frequency truncated terms). We end the article with an appendix gathering results on Sobolev and Besov spaces, the process of diagonalization of System \eqref{systdisp}, and the new Strichartz estimates that allow us to reach this level of precision.

\section{Proof of the results}

\subsection{The limit system}

If $\tUoqg$ is as described in Theorem \ref{Th1}, we regularize it by introducing, for $\lambda>0$ (where $\chi$ is the smooth cut-off function introduced in the appendix)
$$
\tUoqg^\lambda \overset{def}{=} \chi(\frac{|D|}{\lambda}) \tUoqg.
$$
Then $\tUoqg^\lambda \in H^1$ and applying Theorem $2$ from \cite{FC2} there exists a unique global solution $\tUqg^\lambda \in \dot{E}^0 \cap \dot{E}^1$ to System \eqref{QG3} and thanks to Lemma $2.1$ from \cite{FC4} we apply \eqref{estimQG} to $\tUqg^\lambda$ and for all $t\in \R_+$ (taking $\C_0=\max(1,\|\tUoqg\|_{H^{\frac12+\delta}})$):
\begin{multline}
 \|\tUqg^\lambda \|_{L_t^\infty H^{\frac12 + \delta}}^2+\min(\nu, \nu')\int_0^t \|\n \tUqg^\lambda (\tau)\|_{H^{\frac12 + \delta}}^2 d\tau\\
 \leq \Cdn \|\chi(\frac{|D|}{\lambda})\tUoqg\|_{H^{\frac12 + \delta}}^2 \max(1, \|\chi(\frac{|D|}{\lambda})\tUoqg\|_{H^{\frac12 + \delta}}^{\frac{1}{\delta}})\\
 \leq \Cdn \|\tUoqg\|_{H^{\frac12 + \delta}}^2 \max(1, \|\tUoqg\|_{H^{\frac12 + \delta}}^{\frac{1}{\delta}}) \leq \Cdn \C_0^{2+\frac1{\delta}}.
  \label{estimQGlamb}
\end{multline}
Then (taking $\lambda=n$) we prove that $(\tUqg^n)_{n\in \N^*}$ is a Cauchy sequence in $E^{\frac12 + \delta}=\dot{E}^0\cap \dot{E}^{\frac12 + \delta}$. For $n\geq m$, let us define $\td= \tUqg^n-\tUqg^m$, which satisfies the following system:
\begin{equation}
 \begin{cases}
  \displaystyle{\d_t \td-\G \td=-\cQ\left(\tUqg^n \cdot \n \td + \td\cdot \n \tUqg^m\right),}\\
  \displaystyle{\tilde{\delta}_{n,m|t=0}=\big(\chi(\frac{|D|}{n})-\chi(\frac{|D|}{m})\big) \tUoqg.}
 \end{cases}
\end{equation}
For any $s\in[0,\frac12 + \delta]$, taking the $\Hs $-innerproduct and then using the classical Sobolev product laws (see Proposition \ref{HProd}), we get ($(s_1,s_2)\in\{(1,s-\frac12), (s,\frac12)\}$):
\begin{multline}
 \frac12 \frac{d}{dt} \|\td\|_{\Hs }^2 + \nu_0 \|\n \td\|_{\Hs }^2 \leq C\|\tUqg^n \cdot \n \td+ \td\cdot \n \tUqg^m\|_{\dot{H}^{s-1}} \|\td\|_{\dot{H}^{s+1}}\\
 \leq C\left(\|\tUqg^n\|_{\dot{H}^1} \|\td\|_{\Hs }^{\frac12} \|\td\|_{\dot{H}^{s+1}}^{\frac32} +\|\tUqg^m\|_{\dot{H}^\frac32} \|\td\|_{\Hs } \|\td\|_{\dot{H}^{s+1}} \right)\\
 \leq \frac{\nu_0}{2} \|\n \td\|_{\Hs }^2 +\frac{C}{\nu_0} \|\td\|_{\Hs }^2 \left(\|\n\tUqg^m\|_{\dot{H}^\frac12}^2 +\frac{1}{\nu_0^2}\|\tUqg^n\|_{\dot{H}^\frac12}^2 \|\n\tUqg^n\|_{\dot{H}^\frac12}^2 \right).
\end{multline}
Thanks to the Gronwall lemma and using \eqref{estimQGlamb}, we obtain that
$$
\|\td\|_{E^{\frac12 +\delta}}^2 \leq \|\td(0)\|_{H^{\frac12 +\delta}}^2 e^{\frac{\Cdn}{\nu_0^2} \C_0^{2+\frac{1}{\delta}}\left(1+\frac{\Cdn}{\nu_0^2} \C_0^{2+\frac{1}{\delta}}\right)}.
$$
As $\|\td(0)\|_{H^{\frac12 +\delta}}$ goes to zero when $m=\min(n,m)$ goes to infinity, the sequence is Cauchy and if we denote $\tUqg$ its limit in $E^{\frac12 +\delta}$, we immediately get that it solves System \eqref{QG3} and satisfies the expected estimates. $\blacksquare$ 
\\

As an immediate consequence we easily bound $G^{b,l}$ (introduced with the auxiliary systems) as follows:
\begin{prop} \label{estimGlb}
 \sl{There exists a constant $C_F>0$ such that for all $\delta \in]0,\frac12]$ and $s\in[0,\frac12+\delta]$,
 \begin{equation}
 \begin{cases}
 \vspace{2mm}
  \displaystyle{\int_0^\infty \|G^b(\tau)\|_{\Hs } d\tau \leq \frac{C_F}{\nu_0} \Cdn \C_0^{2+\frac{1}{\delta}},}\\
  \displaystyle{\int_0^\infty \|G^l(\tau)\|_{\dot{H}^{s-1}}^2 d\tau \leq C_F \frac{|\nu-\nu'|^2}{\nu_0} \Cdn \C_0^{2+\frac{1}{\delta}}.}
 \end{cases}
 \end{equation}
}
\end{prop}
\begin{rem}
 \sl{In \cite{FC2} the previous terms were estimated for any $s\in[0,1]$ with $\|\tUoqg\|_{H^1}$.}
\end{rem}

\textbf{Proof of Proposition \ref{estimGlb} :} $G^l$ is estimated as in \cite{FC2}, and for $G^b$, as we wish to use only $\frac12+\delta$ derivatives on $\tUoqg$, a much better way than in \cite{FC2} is to write (thanks to the Bony decomposition, see appendix for details):
\begin{multline}
 \|G^b\|_{\Hs } \leq C_F \|\tUqg \cdot \nabla \tUqg \|_{\Hs } \leq C_F \|\div(\tUqg \otimes \tUqg)\|_{\Hs }\\
 \leq C_F \left(2\|T_{\tUqg} \tUqg\|_{\dot{H}^{s+1}} +\|R(\tUqg, \tUqg)\|_{\dot{H}^{s+1}}\right)\\
 \leq C_F \left(2 \|\tUqg\|_{L^\infty} +\|\tUqg\|_{\dot{B}_{\infty, \infty}^0}\right) \|\tUqg\|_{\dot{H}^{s+1}}.
\end{multline}
Then using the injection $\dot{B}_{\infty, 1}^0 \hookrightarrow L^\infty$ together with the Bernstein lemma and the following result (whose proof is close to Lemma $5$ from \cite{FCestimLp}):
\begin{lem} \label{majBs21}
 \sl{For any $\aa, \beta>0$ there exists a constant $C_{\aa, \beta}>0$ such that for any $u\in \dot{H}^{s-\aa} \cap \dot{H}^{s+\beta}$, then $u\in\dot{B}_{2,1}^s$ and:
\begin{equation}
 \|u\|_{\dot{B}_{2,1}^s} \leq C_{\aa, \beta} \|u\|_{\dot{H}^{s-\aa}}^{\frac{\beta}{\aa + \beta}} \|u\|_{\dot{H}^{s+\beta}}^{\frac{\aa}{\aa + \beta}}.
\end{equation}
 }
\end{lem}
we obtain that: 
\begin{equation}
 2 \|\tUqg\|_{L^\infty} +\|\tUqg\|_{\dot{B}_{\infty, \infty}^0} \leq 3\|\tUqg\|_{\dot{B}_{2, 1}^{\frac32}} \leq C\|\tUqg\|_{\dot{H}^{\frac 32-\delta}}^{\frac12} \|\tUqg\|_{\dot{H}^{\frac 32+\delta}}^{\frac12},
\end{equation}
and we end up with (using also \eqref{estimQG}):
\begin{equation}
 \int_0^\infty \|G^b\|_{\Hs } d\tau \leq C_F \|\n \tUqg\|_{L^2 \dot{H}^{\frac 12-\delta}}^{\frac12} \|\n \tUqg\|_{L^2 \dot{H}^{\frac 12+\delta}}^{\frac12} \|\n \tUqg\|_{L^2 \dot{H}^{s}} \leq \frac{C_F}{\nu_0}\Cdn \C_0^{2+\frac{1}{\delta}}. \blacksquare
\end{equation}

\subsection{The case $\nu=\nu'$}

\subsubsection{Estimates for $\We$}

Let us first focus on the linear system \eqref{We}. Let us recall that thanks to Proposition \ref{estimGlb} we obtain that (see \cite{FC2} for details) for any $s\in[\frac12, \frac12 +\delta]$,
\begin{equation}
 \|\We\|_{\dot{E}^s}^2 \leq \left(\|\Uoosc\|_{\Hs }^2 + \frac12 \int_0^t \|G^b(\tau)\|_{\Hs }d\tau\right)e^{\frac12 \int_0^t \|G^b(\tau)\|_{\Hs }} \leq \Do \left(\|\Uoosc\|_{\Hs }^2 + 1\right),
\label{estimWe1}
 \end{equation}
with
$$
\Do\overset{def}{=} \frac{C_F}{\nu_0}\Cdn \C_0^{2+\frac1\delta} e^{\frac{C_F}{\nu_0}\Cdn \C_0^{2+\frac1\delta}}.
$$
One of the main ingredients is to provide a generalization of the Strichartz estimates obtained in \cite{FC3}. Our new Strichartz estimates are much more flexible and we refer to the appendix for the most general formulation (see Propositions \ref{Estimdispnu} and \ref{Estimdispnu2}). We also postpone to the end of the next section the precise statement of the Strichartz estimates we will use.

\subsubsection{Energy estimates}

As explained in section \ref{sol}, we already have a local strong solution $\Ue$ whose lifespan will be denoted as $T_\ee^*$. As seen in the previous section $\tUqg$ and $\We$ exist globally, and then $\de$ is well defined in $\dot{E}_T^{\frac12} \cap \dot{E}_T^{\frac12 +\delta}$ for all  $T<T_\ee^*$ and we can perform for any $s\in[\frac12, \frac12+\eta \delta]$ the innerproduct in $\Hs $ of System \eqref{GE} with $\de$. We have to bound each term from the right-hand side.

Let us begin with the easiest terms, namely $F_1$, $F_2$ and $F_3$: thanks to the classical Sobolev product laws ($(s_1,s_2)=(\frac12,s)$, see Proposition \ref{HProd}), we obtain that:
\begin{equation}
 |(F_1|\de)_{\Hs}|\leq \|\de\cdot \n \de\|_{\dot{H}^{s-1}} \|\de\|_{\dot{H}^{s+1}} \leq C \|\de\|_{\dot{H}^{\frac12}} \|\de\|_{\dot{H}^{s+1}}^2,
 \label{En1}
 \end{equation}
Similarly we obtain that
\begin{equation}
 \begin{cases}
  \vspace{0.2cm}
  \displaystyle{|(F_2|\de)_{\Hs}|\leq C\|\n \tUqg\|_{\dot{H}^{\frac12}} \|\de\|_{\Hs} \|\de\|_{\dot{H}^{s+1}} \leq \frac{\nu}{16} \|\de\|_{\dot{H}^{s+1}}^2 +\frac{C}{\nu} \|\n \tUqg\|_{\dot{H}^{\frac12}}^2 \|\de\|_{\Hs}^2,}\\
  \displaystyle{|(F_3|\de)_{\Hs}|\leq C\|\tUqg\|_{\dot{H}^1} \|\de\|_{\dot{H}^{s+\frac12}} \|\de\|_{\dot{H}^{s+1}} \leq \frac{\nu}{16} \|\de\|_{\dot{H}^{s+1}}^2 +\frac{C}{\nu^3} \|\tUqg\|_{\dot{H}^1}^4 \|\de\|_{\Hs}^2.}
 \end{cases}
 \label{En23}
\end{equation}
Compared to \cite{FC2, FC3} we cannot use for the other $F_i$ the same methods which would produce (after using the Gronwall lemma) a coefficient of the form $e^{\|\We\|_{\dot{E}^s}}$ which would ruin our efforts to allow large initial blow up for the oscilating part (which could only be of size $(-\ln \ee)^\beta$). We need to estimate carefully these terms and especially use as much as possible the new Strichartz estimates (giving positive powers of $\ee$ thanks to Proposition \ref{Estimdispnu}) and the least possible basic energy estimates on $\We$ (that produce $\ee^{-\gamma}$ from \eqref{estimWe1}).

The most obvious way would be to use the paraproduct and remainder laws (see appendix). For example with $F_7$, as $s-1<0$, we have:
\begin{multline}
|(F_7|\de)_{\Hs}|\leq \|\We \cdot \nabla \tUqg\|_{\dot{H}^{s-1}} \|\de\|_{\dot{H}^{s+1}}\\
\leq C\left(\|T_{\We} \nabla \tUqg\|_{\dot{H}^{s-1}}+ \|T_{\nabla \tUqg} \We\|_{\dot{H}^{s-1}} +\|\div(R(\We, \tUqg))\|_{\dot{H}^{s-1}}\right) \|\de\|_{\dot{H}^{s+1}}\\
\leq C\left(\|\We\|_{L^\infty} \|\nabla \tUqg\|_{\dot{H}^{s-1}}+ \|\nabla \tUqg\|_{\dot{H}^{s-1}} \|\We\|_{\dot{B}_{\infty, \infty}^0} +\|\We\|_{\dot{B}_{\infty, \infty}^0} \|\tUqg\|_{\dot{H}^s}\right) \|\de\|_{\dot{H}^{s+1}}\\
\leq C \|\We\|_{\dot{B}_{\infty, 1}^0} \|\tUqg\|_{\dot{H}^s} \|\de\|_{\dot{H}^{s+1}} \leq \frac{\nu}{16} \|\de\|_{\dot{H}^{s+1}}^2 +\frac{C}{\nu} \|\We\|_{\dot{B}_{\infty, 1}^0}^2 \|\tUqg\|_{\dot{H}^s}^2.
\end{multline}
This result could be also usable for $F_5$ but to deal with $\|\We\|_{L^p \dot{B}_{\infty, 1}^0}$ from Proposition \ref{Estimdispnu} we would have to use Lemma \ref{majBs21} which would force us to have a slightly smaller range for $\gamma$. More important, for $F_8$ this method would force us to ask $\gamma<\frac{\delta}{4}$, which is clearly not optimal.

Finally, the most important problem is that the previous estimates cannot be used to estimate $F_4$ and $F_6$: indeed for instance if we wish to estimate $F_6$ the same way:
$$
\|F_6\|_{\dot{H}^{s-1}} \leq C\left(\|T_{\tUqg} \nabla \We\|_{\dot{H}^{s-1}}+ \|T_{\nabla \We} \tUqg\|_{\dot{H}^{s-1}} +\|\div(R(\tUqg, \We))\|_{\dot{H}^{s-1}}\right),
$$
and the first paraproduct (see the appendix for the Bony decomposition) leads to an obstruction as the only possibilities to estimate it are (for $\beta>s$):
\begin{equation}
 \|T_{\tUqg} \nabla \We\|_{\dot{H}^{s-1}} \leq C \begin{cases}
                                               \|\tUqg\|_{L^\infty} \|\We\|_{\dot{H}^s},\\
                                               \|\tUqg\|_{\dot{H}^{s-\beta}} \|\We\|_{\dot{B}_{\infty, \infty}^{\beta}},
                                              \end{cases}
                                              \label{pb}
\end{equation}
In the first estimate each term is well defined but the $\dot{H}^s$-norm of $\We$ produces negative powers of $\ee$, and in the second one the first term is not defined ($\tUqg$ is not defined for negative regularities). It is possible to deal with this term using the same idea as in \cite{FC2} ( with $a,b\geq 1$ so that $\frac1{a}+\frac1{b}=1$),
\begin{multline}
 \int_0^t \|\tUqg \cdot \nabla \We\|_{\dot{H}^{s-1}}^2 d\tau \leq C \int_0^t\|\tUqg \cdot \nabla \We\|_{L^2} \|\tUqg \cdot \nabla \We\|_{\dot{H}^{2(s-1)}} d\tau\\
 \leq \|\tUqg\|_{L^\infty L ^2} \|\nabla \We\|_{L^a L^\infty} \|\tUqg\|_{L^b \dot{H}^{s+\frac12}} \|\We\|_{L^\infty \dot{H}^s},
\end{multline}
and due to the gradient pounding on $\We$, the most interesting use of Proposition \ref{Estimdispnu} consists in choosing $a$ as close as possible to $1$, which implies that $b$ is very large. As $s+\frac12 \geq 1$, this forces us to use \eqref{estimQG} for regularity index close to $1$ (in this case it would be necessary to require that $\tUoqg \in H^s$ with $s$ close to 1), which was something we wished to avoid as we only consider indices  $s\leq\frac12 + \delta$. Moreover it would also produce a clearly non-optimal decrease in $\ee$.

Finally both of these two methods fail for $F_4$: the former for the same reason as for $F_6$, and the latter as we cannot consider $\|\de\|_{L^2}$: there is a lack of derivatives pounding on $\de$.

To overcome this lack of derivatives, we will distribute them differently among the whole $\dot{H}^s$-innerproduct. We will do this for all the last five external force terms and the idea will be to do as in \cite{FCestimLp, FCpochesLp} and deal with the non-local operator $|D|^s$ applied to a product and dispatch $s$ derivatives on $\de$ and obtain something close to the second line of \eqref{pb}. More precisely, we directly deal with the innerproduct as follows:
\begin{equation}
 \Big|(F_4\big|\de)_{\Hs}\Big|= \Big|\big(\div(\de \otimes \We)\big|\de\big)_{\Hs}\Big| =\Big|\big(|D|^s(\de \cdot \We)\big| |D|^s \n \de\big)_{L_2}\Big|.
\label{PSs}
 \end{equation}
The nonlocal operator $|D|^s$ can be written as a singular principal value integral (we refer to \cite{Stein, Cordoba, TH2, TH3, FCestimLp, FCpochesLp}) and when the index $s$ lies in $]0,1[$ (which is the case here as $s$ is close to $\frac12$) it is a classical singular integral:
$$
|D|^s f(x)=C_s \int_{\R^3} \frac{f(x)-f(y)}{|x-y|^{3+s}} dy =C_s \int_{\R^3} \frac{f(x)-f(x-y)}{|y|^{3+s}} dy.
$$
Let us recall that an equivalent formulation of the Besov norm involves translations as stated in the following result:
\begin{thm}\sl{(\cite{Dbook}, $2.36$)
 Let $s \in ]0,1[$ and $p,r\in [1,\infty]$. There exists a constant $C$ such that for any $u\in \dot{B}_{p,r}^s$,
$$
C^{-1} \|u\|_{\dot{B}_{p,r}^s}\leq \|\frac{\|u(\cdot -y) -u(\cdot)\|_{L^p}}{|y|^s}\|_{L^r (\R^d; \frac{dy}{|y|^d})} \leq C \|u\|_{\dot{B}_{p,r}^s}.
$$}
\end{thm}
From this we can prove exactly as in \cite{FCpochesLp} (see section $A.3.1$ there) the following result:
\begin{prop}
\sl{For any $s\in]0,1[$ and any smooth functions $f,g$ we can write:
$$
|D|^s(fg)= (|D|^sf)g +f|D|^sg +M_s(f,g),
$$
where the bilinear operator $M_s$ is defined for all $x\in \R^3$ as:
\begin{equation}
M_s(f,g)(x) =\int_{\R^3} \frac{\big(f(x)-f(x-y)\big) \big(g(x)-g(x-y)\big)}{|y|^{3+s}} dy.
\label{defM}
\end{equation}
Moreover there exists a constant $C_s$ such that for all $f,g$ and all $p,p_1, p_2,r_1,r_2\in[1,\infty]$ and $s_1,s_2>0$ satisfying:
$$
\displaystyle{\frac{1}{p}= \frac{1}{p_1} +\frac{1}{p_2}, \quad 1= \frac{1}{r_1} +\frac{1}{r_2}, \quad s_1+s_2=s},
$$
then we have
\begin{equation}
\|M_s(f,g)\|_{L^p} \leq C_s \|f\|_{\displaystyle{\dot{B}_{p_1,r_1}^{s_1}}}\|g\|_{\displaystyle{\dot{B}_{p_2,r_2}^{s_2}}}.
\end{equation}
}
\label{propM2}
\end{prop}
\begin{rem}
 \sl{The additional term $M_s$ allows us to freely dispatch the derivatives as desired provided that $s_1,s_2>0$, which will force us to spend a small extra amount of derivative in order to meet these conditions. So even if it is not possible to use Proposition \ref{propM2} for $(s_1,s_2)=(s,0)$, our method will enable us to do nearly as if we could estimate $\|M_s(\de,\We)\|_{L^2}$ by $\|\de\|_{\dot{H}^\frac12}\||D|^s\We\|_{L^6}$.}
\end{rem}
More precisely for a small $\aA>0$, instead of \eqref{PSs}, we will write (also using the Sobolev injections):
\begin{multline}
 \Big|(F_6\big|\de)_{\Hs}\Big|= \Big|\big(\div(\tUqg \otimes \We)\big|\de\big)_{\Hs}\Big| =\Big|\big(|D|^{s+\aA}(\tUqg \cdot \We)\big| |D|^{s-\aA} \n \de\big)_{L_2}\Big|\\
 \leq \|(|D|^{s+\aA}\tUqg)\cdot\We +\tUqg\cdot|D|^{s+\aA}\We +M_{s+\aA}(\tUqg,\We)\|_{L^\frac{6}{3+2\aA}} \cdot \||D|^{s-\aA} \n \de\|_{L^\frac{6}{3-2\aA}}\\
 \leq C\left(\||D|^{s+\aA}\tUqg\|_{L^2} \|\We\|_{L^\frac{3}{\aA}} +\|\tUqg\|_{L^3}\||D|^{s+\aA}\We\|_{L^\frac{6}{1+2\aA}} +\|\tUqg\|_{\dot{H}^s}\|\We\|_{\dot{B}_{\frac{3}{\aA},2}^{\aA}} \right)\\
 \times \||D|^{s-\aA}\n \de\|_{\dot{H}^{\aA}}\\
 \leq C\left(\|\tUqg\|_{\dot{H}^{s+\aA}} \|\We\|_{L^\frac{3}{\aA}} +\|\tUqg\|_{\dot{H}^\frac12} \||D|^{s+\aA}\We\|_{L^\frac{6}{1+2\aA}} +\|\tUqg\|_{\dot{H}^s} \|\We\|_{\dot{B}_{\frac{3}{\aA},2}^{\aA}}\right) \cdot \|\de\|_{\dot{H}^{s+1}}\\
 \leq \frac{\nu}{16} \|\de\|_{\dot{H}^{s+1}}^2 +\frac{C}{\nu} \Bigg( \|\tUqg\|_{\dot{H}^s}^{2(1-\aA)} \|\tUqg\|_{\dot{H}^{s+1}}^{2\aA} \|\We\|_{L^\frac{3}{\aA}}^2 +\|\tUqg\|_{\dot{H}^\frac12}^2 \||D|^{s+\aA}\We\|_{L^\frac{6}{1+2\aA}}^2\\
 +\|\tUqg\|_{\dot{H}^s}^2 \|\We\|_{\dot{B}_{\frac{3}{\aA},2}^{\aA}}^2\Bigg).
 \label{En67}
 \end{multline}
\begin{rem}
 \sl{Notice that as $\de, \We, \tUqg$ are divergence-free, we will systematically (thanks to integration by parts) transfer the divergence as a gradient on the right-hand part of the innerproduct, and as a consequence the computations are the same respectively for $F_4$ and $F_5$, and for $F_6$ and $F_7$.}
\end{rem}
Let us continue with $F_4$, by the classical Sobolev interpolation and Young estimates, we can write that (for $\aB>0$ small):
\begin{multline}
 \Big|(F_4\big|\de)_{\Hs}\Big|= \Big|\big(\div(\de \otimes \We)\big|\de\big)_{\Hs}\Big| =\Big|\big(|D|^{s+\aB}(\de \cdot \We)\big| |D|^{s-\aB} \n \de\big)_{L_2}\Big|\\
 \leq C\|\left(|D|^{s+\aB}\de)\cdot\We +\de\cdot|D|^{s+\aB}\We +M_{s+\aB}(\de,\We\right)\|_{L^\frac{6}{3+2\aB}} \cdot \||D|^{s-\aB} \n \de\|_{L^\frac{6}{3-2\aB}}\\
 \leq C\left(\||D|^{s+\aB}\de\|_{L^2} \|\We\|_{L^\frac{3}{\aB}} +\|\de\|_{L^3}\||D|^{s+\aB}\We\|_{L^\frac{6}{1+2\aB}} +\|\de\|_{\dot{H}^s}\|\We\|_{\dot{B}_{\frac{3}{\aB},2}^{\aB}} \right) \cdot \||D|^{s-\aB} \n \de\|_{\dot{H}^{\aB}}\\
  \leq C \|\de\|_{\dot{H}^s}^{1-\aB} \|\de\|_{\dot{H}^{s+1}}^{1+\aB} \|\We\|_{L^\frac{3}{\aB}} + C\left( \|\de\|_{L^3}\||D|^{s+\aB} \We\|_{L^\frac{6}{1+2\aB}} +\|\de\|_{\dot{H}^s}\|\We\|_{\dot{B}_{\frac{3}{\aB},2}^{\aB}}\right) \|\de\|_{\dot{H}^{s+1}}\\
  \leq \frac{\nu}{16} \|\de\|_{\dot{H}^{s+1}}^2 +C \|\de\|_{\dot{H}^s}^2 \left(\frac{1}{\nu^{\frac{1+\aB}{1-\aB}}}\|\We\|_{L^\frac{3}{\aB}}^\frac{2}{1-\aB} +\frac{1}{\nu}\|\We\|_{\dot{B}_{\frac{3}{\aB},2}^{\aB}}^2 \right)\\
  +\frac{C}{\nu}\|\de\|_{\dot{H}^\frac12}^2 \||D|^{s+\aB}\We\|_{L^\frac{6}{1+2\aB}}^2.
 \label{En45}
\end{multline}
Finally we estimate $F_8$ with the same method, but the term $M_{s+\aC}(\We,\We)$ has to be estimated differently (otherwise we end up with the same problem as explained in the beginning of this section): instead of estimating it as for the other terms by $\|\We\|_{\dot{H}^s} \|\We\|_{\dot{B}_{\frac{3}{\aC},2}^{\aC}}$ (the first term being $L^\infty$, and the second $L^2$ in time), we will estimate it by
$$
\|\We\|_{\dot{H}^{s+\aC-\bd}}\|\We\|_{\dot{B}_{\frac{3}{\aC},2}^{\bd}},
$$
for small enough $\aC,\beta>0$ so that the first term keeps $L^\infty$ in time and the second one is $L^2$ (we try to be as close as possible to the forbidden choice $\beta=0$). As we will precise below, dealing with $\|\We\|_{L^\infty \dot{H}^{s}}^{2(1-\aC)} \|\We\|_{L^2 \dot{H}^{s}}^{2\aC} \|\We\|_{L^2 L^\frac{3}{\aC}}^2$ (for the first term) will only lead to $\gamma<\frac{\delta}4$, whereas $\|\We\|_{L^\infty \dot{H}^{s+\aC}}^2 \|\We\|_{L^2 L^\frac{3}{\aC}}^2$  will allow us to reach $\gamma<\frac{\delta}2$. For the same reason we will estimate the other term by $\|\We\|_{L^2 \dot{B}_{\frac{3}{\aC},2}^{\bd}}$ instead of $\|\We\|_{L^\frac{2}{1-\aC} \dot{B}_{\frac{3}{\aC},2}^{\bd}}$. Altough this choice seems very close to the other, it allows us to use a smaller $p$ in the Strichartz estimates, which allows a slightly wider range for $\theta$ helping us to reach $\gamma<\frac{\delta}{2}$ instead of $\gamma<\frac{\delta}{4}$. Once more, we try to obtain as close as possible to what we would get if it Proposition \ref{propM2} could be applied for $s_1=s+\aC$ and $s_2=0$.
\begin{multline}
 \Big|(F_8\big|\de)_{\Hs}\Big| \leq\||D|^{s+\aC}(\We \otimes \We)\|_{L^\frac{6}{3+2\aC}} \||D|^{s-\aC} \n \de\|_{L^\frac{6}{3-2\aC}}\\
 \leq \left(2\||D|^{s+\aC}\We\|_{L^2} \|\We\|_{L^\frac{3}{\aC}} +\|\We\|_{\dot{H}^{s+\aC-\bd}}\|\We\|_{\dot{B}_{\frac{3}{\aC},2}^{\bd}} \right) \cdot \||D|^{s-\aC}\n\de\|_{\dot{H}^{\aC}}\\
 \leq \frac{\nu}{16} \|\de\|_{\dot{H}^{s+1}}^2 +\frac{C}{\nu} \left( \|\We\|_{\dot{H}^{s+\aC}}^2 \|\We\|_{L^\frac{3}{\aC}}^2 +\|\We\|_{\dot{H}^{s+\aC-\bd}}^2 \|\We\|_{\dot{B}_{\frac{3}{\aC},2}^{\bd}}^2 \right).
 \label{En8}
 \end{multline}
We can now gather all the external force term estimates \eqref{En1}, \eqref{En23}, \eqref{En45}, \eqref{En67}, \eqref{En8} and taking the $\dot{H}^s$-innerproduct of System \eqref{GE} with $\de$, we obtain that for all $s\in[\frac12, \frac12 +\eta\delta]$ and all $t<T_\ee^*$:
\begin{multline}
 \frac12 \frac{d}{dt} \|\de\|_{\dot{H}^s}^2 +\nu \|\n \de\|_{\dot{H}^s}^2 \leq \left(C\|\de\|_{\dot{H}^{\frac12}} +8\frac{\nu}{16}\right) \|\n \de\|_{\dot{H}^s}^2\\
 +\frac{C}{\nu} \|\de\|_{\dot{H}^s}^2 \Bigg\{ \|\n \tUqg\|_{\dot{H}^{\frac12}}^2 (1+\frac{1}{\nu^2}\|\tUqg\|_{\dot{H}^{\frac12}}^2) +\frac{1}{\nu^{\frac{2\aB}{1-\aB}}}\|\We\|_{L^\frac{3}{\aB}}^\frac{2}{1-\aB} +\|\We\|_{\dot{B}_{\frac{3}{\aB},2}^{\aB}}^2 \Bigg\}\\
 +\frac{C}{\nu} \Bigg[ \|\tUqg\|_{\dot{H}^s}^{2(1-\aA)} \|\tUqg\|_{\dot{H}^{s+1}}^{2\aA} \|\We\|_{L^\frac{3}{\aA}}^2 +\|\tUqg\|_{\dot{H}^\frac12}^2 \||D|^{s+\aA}\We\|_{L^\frac{6}{1+2\aA}}^2 +\|\tUqg\|_{\dot{H}^s}^2 \|\We\|_{\dot{B}_{\frac{3}{\aA},2}^{\aA}}^2\\
 +\|\de\|_{\dot{H}^\frac12}^2 \||D|^{s+\aB}\We\|_{L^\frac{6}{1+2\aB}}^2 +\|\We\|_{\dot{H}^{s+\aC}}^2 \|\We\|_{L^\frac{3}{\aC}}^2 +\|\We\|_{\dot{H}^{s+\aC-\bd}}^2 \|\We\|_{\dot{B}_{\frac{3}{\aC},2}^{\bd}}^2 \Bigg].
\end{multline}
In order to perform the bootstrap argument (we refer to in \cite{FC2, FC3}), let us now define
\begin{equation}
 T_\ee \overset{def}{=} \sup \{t\in[0,T_\ee^*[, \quad \forall t'\leq t, \|\de(t')\|_{\dot{H}^\frac12} \leq \frac{\nu}{4C}\}.
\end{equation}
Due to the assumptions, $\|\de(0)\|_{H^{\frac12+\delta}}\leq \Co \ee^{\aa_0}$ so that we are sure that $T_\ee>0$ if $\ee \leq \left(\frac{\nu}{8C\Co}\right)^\frac{1}{\aa_0}$. Thanks to the Gronwall and Young estimates, and estimating the first terms in the last block as follows:
\begin{multline}
 \int_0^\infty \|\tUqg\|_{\dot{H}^s}^{2(1-\aA)} \|\tUqg\|_{\dot{H}^{s+1}}^{2\aA} \|\We\|_{L^\frac{3}{\aA}}^2 d\tau\\
 \leq \left(\int_0^\infty \|\tUqg\|_{\dot{H}^{s+1}}^2 d\tau \right)^{\aA} \left(\int_0^\infty \|\We\|_{L^\frac{3}{\aA}}^\frac{2}{1-\aA} \|\tUqg\|_{\dot{H}^s}^2 d\tau \right)^{1-\aA},
\end{multline}
we can now state that for all $s\in[\frac12, \frac12 +\eta\delta]$ and all $t\leq T_\ee$, we have (as $\We$ and $\tUqg$ are globally defined, each time integral in the right-hand side is over $\R_+$):
\begin{multline}
 \|\de(t)\|_{\dot{H}^s}^2 +\frac{\nu}{2} \int_0^t \|\n \de(\tau)\|_{\dot{H}^s}^2 d\tau \leq \Bigg[ \|\de(0)\|_{\dot{H}^s}^2 +\frac{C}{\nu} \Bigg(\|\tUqg\|_{L^\infty\dot{H}^s}^{2(1-\aA)} \|\tUqg\|_{L^2 \dot{H}^{s+1}}^{2\aA} \|\We\|_{L^\frac{2}{1-\aA} L^\frac{3}{\aA}}^2\\
 +\|\tUqg\|_{L^\infty \dot{H}^\frac12}^2 \||D|^{s+\aA}\We\|_{L^2 L^\frac{6}{1+2\aA}}^2 +\|\tUqg\|_{L^\infty \dot{H}^s}^2 \|\We\|_{L^2 \dot{B}_{\frac{3}{\aA},2}^{\aA}}^2\\
 +\||D|^{s+\aB}\We\|_{L^2 L^\frac{6}{1+2\aB}}^2 +\|\We\|_{L^\infty \dot{H}^{s+\aC}}^2 \|\We\|_{L^2 L^\frac{3}{\aC}}^2 +\|\We\|_{L^\infty \dot{H}^{s+\aC-\bd}}^2 \|\We\|_{L^2 \dot{B}_{\frac{3}{\aC},2}^{\bd}}^2 \Bigg) \Bigg]\\
\times \exp \frac{C}{\nu}\Bigg\{ \|\n \tUqg\|_{L^2 \dot{H}^{\frac12}}^2 (1+\frac{1}{\nu^2}\|\tUqg\|_{L^\infty \dot{H}^{\frac12}}^2) +\frac{1}{\nu^{\frac{2\aB}{1-\aB}}}\|\We\|_{L^\frac{2}{1-\aB} L^\frac{3}{\aB}}^\frac{2}{1-\aB} +\|\We\|_{L^2 \dot{B}_{\frac{3}{\aB},2}^{\aB}}^2 \Bigg\}.
\label{estimCas1}
\end{multline}
It is now about to properly use the new Strichartz estimates we proved in the present article (see the appendix for Proposition \ref{Estimdispnu} and its proof).
\\

Let us begin with the case $(d,p,r,q)=(s+\aa,2,\frac{6}{1+2\aa},2)$, for all $\theta\in ]0,\frac{1-\aa}{1-4\aa}[\cap]0,1]=]0,1]$. Thanks to Proposition \ref{injectionLr} (for more simplicity we will not track the dependency in $\nu$),
\begin{multline}
  \||D|^{s+\aa} \We\|_{L_t^2 L^\frac{6}{1+2\aa}} \leq C \||D|^{s+\aa} \We\|_{\tilde{L}_t^2\dot{B}_{\frac{6}{1+2\aa},2}^0}\\
  \leq C_{F,\nu, p,\theta,\aa} \ee^{\frac{\theta}{12}(1-4\aa)} \left( \|\Uoosc\|_{\dot{H}^{s+\frac{\theta}{6} (1-4\aa)}} +\int_0^t  \|G^b(\tau)\|_{\dot{H}^{s+\frac{\theta}{6} (1-4\aa)}} d\tau \right),
 \end{multline}
and if we choose $\aa\in]0,\frac14[$, and $\theta=\frac{6(\delta+\frac12-s)}{1-4\aa}$ (which is in $]0,1]$ if $\delta \leq s-\frac13-\frac{2\aa}{3}$, recall that $s\sim \frac12$), then we obtain (thanks to Proposition \ref{estimGlb}):
\begin{multline}
  \||D|^{s+\aa} \We\|_{L_t^2 L^\frac{6}{1+2\aa}} \leq C \||D|^{s+\aa} \We\|_{\tilde{L}_t^2\dot{B}_{\frac{6}{1+2\aa},2}^0}\\
  \leq C_{F, \nu, s,\delta,\aa} \ee^{\frac12 (\delta+ \frac12-s)} \left( \|\Uoosc\|_{\dot{H}^{\frac12+\delta}} +\int_0^t  \|G^b(\tau)\|_{\dot{H}^{\frac12+\delta}} d\tau \right)\\
  \leq C_{F, \nu, s,\delta,\aa} \ee^{\frac12 (\delta+ \frac12-s)} \Do(\|\Uoosc\|_{\dot{H}^{\frac12+\delta}}+ 1).
  \label{StriA}
 \end{multline}
Let us continue with the case $(d,p,r,q)=(\aa,2,\frac{3}{\aa},2)$, for all $\theta\in ]0,\frac{\frac12-\frac{\aa}{3}}{1-\frac{4\aa}{3}}[$, if we assume $\aa\in]0,\frac34[$, and choose $\theta=\frac{6\delta}{3-4\aa}$,
\begin{equation}
  \|\We\|_{\tilde{L}_t^2\dot{B}_{\frac{3}{\aa},2}^{\aa}} \leq C_{F, \nu,\delta, \aa} \ee^{\frac{\delta}{2}} \Do \big(\|\Uoosc\|_{\dot{H}^{\frac12+\delta}}+ 1\big).
  \label{StriB}
 \end{equation}
For the case $(d,p,r,q)=(0,\frac{2}{1-\aa},\frac{3}{\aa},2)$, for all $\theta\in ]0,\frac{\frac12-\frac{\aa}{3}}{1-\frac{4\aa}{3}}[$, if $\aa\in]0,\frac34[$, and if we choose $\theta=\frac{6\delta}{3-4\aa}$,
\begin{equation}
  \|\We\|_{L_t^\frac{2}{1-\aa} L^\frac{3}{\aa}} \leq C\|\We\|_{\tilde{L}_t^\frac{2}{1-\aa}\dot{B}_{\frac{3}{\aa},2}^0} \leq C_{F, \nu,\delta, \aa} \ee^{\frac{\delta}{2}} \Do(\|\Uoosc\|_{\dot{H}^{\frac12+\delta}}+ 1).
  \label{StriC}
 \end{equation}
All these estimates are verified for $\aA=\aB=\aa= \frac1{16}$ if $\delta \leq \frac18$. Then we turn to the last two terms from \eqref{En8}, let us begin by the first one: as announced, due to the first factor (estimated thanks to \eqref{estimWe1}), doing as before will only allow us to get $\ee^\frac{\delta}{2} \Do(\|\Uoosc\|_{\dot{H}^{\frac12+\delta}}+ 1)^2$, which leads to $\gamma< \frac{\delta}{4}$. In order to reach the announced bound $\frac{\delta}{2}$, we will try to take a slightly smaller $p$ which will allow us to widen the range for $\theta$. But taking $p=2$ instead of $\frac{2}{1-\aa}$ requires that $\|\We\|_{\dot{H}^{s+\aC}}$ is $L^\infty$, that is we need that $s+\aC\leq \frac12 +\delta$. More precisely with $(d,p,r,q)=(0,2,\frac{3}{\aa},2)$, we have
\begin{equation}
  \|\We\|_{L_t^2 L^\frac{3}{\aC}} \leq C\|\We\|_{\tilde{L}_t^2 \dot{B}_{\frac{3}{\aC},2}^0} \leq C_{F, \nu, \theta,s} \ee^{\frac{\theta}{4}(1-\frac{4\aC}{3})} \Do \big(\|\Uoosc\|_{\dot{H}^{\frac12-\aC+\frac{\theta}{2}(1-\frac{4\aC}{3})}}+ 1\big).
 \end{equation}
and as we want
$$
 \aC+s =\frac12 +\delta =\frac12-\aC+\frac{\theta}{2}(1-\frac{4\aC}{3})
$$
we choose
$$
 (\aC, \theta)=(\delta +\frac12-s, \frac{2(\delta +\aC)}{1-\frac{4\aC}{3}}),
$$
which is possible (according to the condition from Proposition \ref{Estimdispnu}) when $\theta <\frac{\frac12-\frac{\aC}{3}}{1-\frac{4\aC}{3}}$, that is if
\begin{equation}
 \delta <\frac{7s-2}{13},
 \label{Condelta1}
\end{equation}
which is realized (recall that $s\in[\frac12, \frac12 +\eta \delta]$) when $\delta \leq \frac1{10}<\frac{3}{26}$, then we have
\begin{equation}
  \|\We\|_{L_t^2 L^\frac{3}{\aC}} \leq C\|\We\|_{\tilde{L}_t^2 \dot{B}_{\frac{3}{\aC},2}^0} \leq C_{F, \nu, \delta,s} \ee^{\frac12 \Big(2\delta +\frac12 -s\Big)} \Do(\|\Uoosc\|_{\dot{H}^{\frac12+\delta}}+ 1).
  \label{StriD}
 \end{equation}
Now, for the last term, $\aC$ is fixed and we will adjust $\theta$ and $\beta$. For $(d,p,r,q)=(\bd,2,\frac{3}{\aC},2)$, we choose $\theta$ so that the corresponding $\sigma$ (see Proposition \ref{Estimdispnu}) is equal to $\frac12 + \delta$, that is
$$
\frac{\theta}{2}(1-\frac{4\aC}{3})=(2-\beta) \delta+\frac12-s,
$$
which is possible when $\theta\in ]0,\frac{\frac12-\frac{\aC}{3}}{1-\frac{4\aC}{3}}[$ that is $\delta < \frac{7s-2}{13-6\beta}$, which is realized when \eqref{Condelta1} is true (when $\beta \in ]0,1[$). In this case, we end up with
\begin{equation}
  \|\We\|_{\tilde{L}_t^2 \dot{B}_{\frac{3}{\aC},2}^{\bd}} \leq C_{F,\nu,\aa,\delta,s} \ee^{\frac12 \Big((2-\beta)\delta +\frac12 -s\Big)} \Do(\|\Uoosc\|_{\dot{H}^{\frac12+\delta}}+ 1).
  \label{StriE}
 \end{equation}
 Combining \eqref{estimCas1} with all these Strichartz estimates, namely \eqref{StriA}, \eqref{StriB}, \eqref{StriC}, \eqref{StriD} and \eqref{StriE}, we end-up for all $s\in[\frac12, \frac12 +\eta\delta]$, all $\beta>0$ small and all $t\leq T_\ee$ with
 \begin{multline}
 \|\de(t)\|_{\dot{H}^s}^2 +\frac{\nu}{2} \int_0^t \|\n \de(\tau)\|_{\dot{H}^s}^2 d\tau\\
 \leq \Bigg[ \|\de(0)\|_{\dot{H}^s}^2  +\Do\Big( (\ee^{\delta+\frac12-s}+\ee^\delta)(\|\Uoosc\|_{\dot{H}^{\frac12+\delta}}+ 1)^2
 +(\ee^{2\delta+\frac12-s}+\ee^{(2-\beta)\delta+\frac12-s})(\|\Uoosc\|_{\dot{H}^{\frac12+\delta}}+ 1)^4\Big)\Bigg]\\
 \times \exp \Bigg\{\Do\left(1+\ee^\delta (\|\Uoosc\|_{\dot{H}^{\frac12+\delta}}+ 1)^2 +\Big(\ee^\delta (\|\Uoosc\|_{\dot{H}^{\frac12+\delta}}^2+1)\Big)^\frac2{1-\aB}\right) \Bigg\}\\
 \leq \Do \left[\ee^{2\aa_0} + \Big( (\ee^{\delta+\frac12-s}+\ee^\delta)\|\Uoosc\|_{\dot{H}^{\frac12+\delta}}^2 +(\ee^{2\delta+\frac12-s}+\ee^{(2-\beta)\delta+\frac12-s})\|\Uoosc\|_{\dot{H}^{\frac12+\delta}}^4\right]\\
 \times \exp \Bigg\{\Do\Big(1+\ee^\delta \|\Uoosc\|_{\dot{H}^{\frac12+\delta}}^2 +(\ee^\delta \|\Uoosc\|_{\dot{H}^{\frac12+\delta}}^2)^\frac2{1-\aB} \Big)\Bigg\}.
\label{estimCas1b}
\end{multline}
As $s\in[\frac12, \frac12+\eta\delta]$, we can write that:
\begin{multline}
\|\de(t)\|_{\dot{H}^s}^2 +\frac{\nu}{2} \int_0^t \|\n \de(\tau)\|_{\dot{H}^s}^2 d\tau\\
\leq \Do \left[\ee^{2\aa_0} +\ee^{(1-\eta)\delta-2\gamma} +\ee^{(2-\eta)\delta-4\gamma} +\ee^{(2-\eta-\beta)\delta-4\gamma}\right] e^{\Do\Big(1+\ee^{\delta-2\gamma}\Big)},
\end{multline}
so that we need
$$
\gamma <\min\Big((1-\eta)\frac{\delta}{2}, (1-\frac{\eta}{2})\frac{\delta}{2}, (1-\frac{\beta+\eta}{2})\frac{\delta}{2}\Big).
$$
If we fix $\beta=\eta$, the condition is reduced to $\gamma<(1-\eta)\frac{\delta}{2}$, so that if $0<\gamma<\frac{\delta}{2}$, we define $\eta=\frac12 (1-\frac{2\gamma}{\delta})$ (or equivalently $\gamma=(1-2\eta)\frac{\delta}{2}$), then with $\beta=\eta$, for all $s\in [\frac12, \frac12+\eta\delta]$ and $t\leq T_\ee$, we end up with (as soon as $\ee\leq 1$):
\begin{equation}
\|\de(t)\|_{\dot{H}^s}^2 +\frac{\nu}{2} \int_0^t \|\n \de(\tau)\|_{\dot{H}^s}^2 d\tau \leq \Do e^{2\Do} \ee^{2\min(\aa_0, \frac{\eta \delta}{2})}.
\end{equation}
We can now conclude the bootstrap argument: there exists $\ee_0>0$ such that for any $0<\ee<\ee_0$ the previous quantity is bounded by $\left(\frac{\nu}{8C}\right)^2$, so that (in particular for $s=\frac12$) if we assume by contradiction that $T_\ee<T_\ee^*$, then $\|\de\|_{L_{T_\ee}^\infty \dot{H}^\frac12} \leq \frac{\nu}{8C}$, which contradicts the maximality of $T_\ee$ (in this case, we would have $\|\de(T_\ee)\|_{\dot{H}^\frac12} =\frac{\nu}{4C}$). Then $T_\ee=T_\ee^*$ and the previous estimates hold true for any $t<T_\ee^*$, so that by the blowup criterion $T_\ee=T_\ee^*=\infty$ and the previous estimate is true for all $t\geq 0$ and all $s\in[\frac12, \frac12 +\eta \delta]$:
$$
\|\de\|_{\dot{E}^s}\leq \Bo \ee^{\min(\aa_0, \frac{\eta \delta}{2})}.
$$
Finally, to prove the last part of the theorem, we only have to remark that the previous argument is then true for any $s\in[\frac12-\eta \delta, \frac12+\eta \delta]$ when we ask $\delta <\frac{3}{26+14\eta}$ (instead of $\delta <\frac{3}{26}$, see \eqref{Condelta1}), and use Lemma \ref{majBs21}:
$$
\|\de\|_{L^2 L^\infty}\leq \|\de\|_{L^2 \dot{B}_{\infty,1}^0} \leq (\|\de\|_{L^2 \dot{H}^{\frac32-\eta \delta}} \|\de\|_{L^2 \dot{H}^{\frac32+\eta \delta}})^\frac12 \leq \Bo \ee^{\min(\aa_0, \frac{\eta \delta}{2})}.
$$
For $(d,p,r,q)=(0,2,\infty,1)$ and for all $\theta\in]0,\frac12[$, from Proposition \ref{Estimdispnu}:
$$
\|\We\|_{L^2 L^\infty} \leq \Co \ee^\frac{\theta}{4} \left(\|\Uoosc\|_{\dot{B}_{2,1}^{\frac12 +\frac{\theta}{2}}}+\int_0^\infty \|G^b(\tau)\|_{\dot{B}_{2,1}^{\frac12 +\frac{\theta}{2}}} d\tau\right).
$$
Using Lemma \ref{majBs21} with $(\aa,\beta)=(\frac{\theta}{2}, k \frac{\theta}{2})$, and if $\theta=\frac{2\delta}{1+k}$ (for some small $k>0$),
\begin{equation}
 \|\Uoosc\|_{\dot{B}_{2,1}^{\frac12 +\frac{\theta}{2}}}\leq \|\Uoosc\|_{\dot{H}^\frac12}^\frac{k}{1+k} \|\Uoosc\|_{\dot{H}^{\frac12 +(1+k)\frac{\theta}{2}}}^\frac{1}{1+k} \leq \Co \ee^{-\gamma}.
\end{equation}
Choosing $k=\frac{\eta}{1-\eta}$, we get
$$
\|\We\|_{L^2 L^\infty} \leq \Do \ee^{\frac{\delta}{2}(\frac{1}{1+k}-(1-2\eta))}=\Do \ee^\frac{\eta \delta}{2},
$$
and the conclusion follows from the fact that $\Ue-\tUqg=\de +\We$. $\blacksquare$
\begin{rem}
 \sl{Going back to \eqref{estimCas1b}, in the case $s=\frac12$ if we only seek for global well posedness, we retrieve here the same condition as in \cite{IMT}, except for the last term because Proposition \ref{propM2} imposes $\beta>0$, so that the condition for global wellposedness is still $\gamma<\frac{\delta}{2}$. If $\beta$ could reach zero, the conditions would be:
 $$
 \begin{cases}
 \vspace{0.2cm}
  \|\Uoosc\|_{\dot{H}^\frac12\cap \dot{H}^{\frac12+\delta}}^2 \ee^{\delta}\leq c, \mbox{ with $c$ some fixed small constant, if we want global well-posedness,}\\
  \vspace{0.2cm}
  \|\Uoosc\|_{\dot{H}^\frac12\cap \dot{H}^{\frac12+\delta}}^2 \ee^{\delta}\underset{\ee\rightarrow 0}{\rightarrow} 0, \mbox{ if we want that } \de\mbox{ goes to zero,}\\
  \gamma<\frac{\delta}{2}, \mbox{ if we want that } \de \mbox{ goes to zero as a positive power of }\ee \mbox{ (which is what we originally searched for).}
 \end{cases}
 $$
 In our case, due to this $\beta>0$ these three conditions coincide.}
 \label{RqIMT}
\end{rem}

\subsection{The general case}

\subsubsection{Estimates on $\Wet$}

Let us begin by recalling the energy estimates for $\Wet$ (we refer to Theorem \ref{Th3} for $m,M$).
\begin{prop}\label{estimWe}
 \sl{Assume $M<\frac12$, there exist $\ee_0 = \ee_0(\nu, \nu', M)>0$ and $\Bo=\Bo(\Co, \nu, \nu', F,\delta)\geq 1$ such that for any $0<\ee\leq \ee_0$ and $s\in \intervb$, we have:
\begin{equation}
  \|\Wet\|_{L^\infty (\R_+, \Hs )}^2 +\nu_0 \|\Wet\|_{L^2 (\R_+, \dot{H}^{s+1})}^2 \leq \Bo \left(\ee^{-2\gamma} + 1\right).
  \label{estimWes}
 \end{equation}
 }
\end{prop}

\textbf{Proof :} we know from \cite{FC2} that there exists a constant $C_F>0$ such that for any $s\in[0,1]$ and $t\in \R_+$, we have:
 \begin{multline}
  \|\Wet\|_{\dot{E}_t^s} \leq e^{\int_0^t \|G^b(\tau)\|_{\Hs } d\tau} \\
  \times \left(\|\Wet(0)\|_{\Hs }^2 +C_F (1+\ee \Re^2 |\nu-\nu'|)^2 \int_0^t (\|G^b(\tau)\|_{\Hs } +\frac{1}{\nu_0} \|G^l(\tau)\|_{\dot{H}^{s-1}}) d\tau\right).
 \end{multline}
Combined with \eqref{estimvp3+4}, Proposition \ref{estimGlb} allows us to obtain that when $s\in\intervb$:
\begin{multline}
  \|\Wet\|_{\dot{E}_t^s} \leq C_F (1+|\nu-\nu'| \ee \Re^2)^2 e^{\frac{C_F}{\nu_0} \Cdn \C_0^{2+\frac1{\delta}}} \\
  \times \left(\|\cPrr \Uoosc \|_{\Hs }^2 +(\frac{1}{\nu_0} +\frac{|\nu-\nu'|^2}{\nu_0^2}) \Cdn \C_0^{2+\frac1{\delta}} \right).
 \end{multline}
We have $|\nu-\nu'| \ee \Re^2 \leq 1$ as soon as $M>\frac12$ and $\ee \leq \ee_0 =|\nu-\nu'|^{\frac{-1}{1-2M}}$ which leads to \eqref{estimWes}. %choosing
%$$
%\Bo \geq 4C_F e^{\frac{C_F \Cdn \Co^{2+\frac1{\delta}}}{\nu_0}} (1+\frac{1}{\nu_0} +\frac{|\nu-\nu'|^2}{\nu_0^2}) \Cdn \Co^{2+\frac1{\delta}}. \blacksquare
%$$

\subsubsection{Estimates on $\de$}

As explained in the previous section (see also \cite{FC2, FCpochesLp}), as $\Uoe\in \Hs $ for all $s\in \intervb$, in particular it lies in $\dot{H}^\frac12$ and thanks to the Fujita and Kato theorem there exists a unique local strong solution $\Ue\in L_T^\infty \dot{H}^\frac12 \cap L_T^2 \dot{H}^\frac32$ for all $0<T<T_\ee^*$ where $T_\ee^*>0$ denotes the maximal lifespan. In addition, if $T_\ee^*$ is finite then we have:
$$
\int_0^{T_\ee^*} \|\n \Ue(\tau)\|_{\dot{H}^\frac{1}{2}(\R^3)}^2 d\tau=\infty.
$$
Moreover, as our initial data enjoys additionnal regularity properties, they are transmitted to the solution: for all $s\in \intervb$ and $T<T_\ee^*$, $\Ue \in L_T^\infty \Hs  \cap L_T^2 \dot{H}^{s+1}$.
As before, with a view to a bootstrap argument, let us now define
\begin{equation}
 T_\ee \overset{def}{=} \sup \{t\in[0,T_\ee^*[, \quad \forall t'\leq t, \|\de(t')\|_{\dot{H}^\frac12} \leq \frac{\nu}{4C}\}.
\end{equation}
Thanks to \eqref{estimfext1}, we are sure that $\|\de(0)\|_{\dot{H}^\frac12}\leq \frac{\nu}{8C}$ (and then $T_\ee>0$) if $\ee$ is small enough. Assuming that $T_\ee<T_\ee^*$, the computations from the previous case imply that, for all $s\in \interv$, and all $t\leq T_\ee$,
\begin{multline}
 \|\de(t)\|_{\dot{H}^s}^2 +\frac{\nu_0}{2} \int_0^t \|\n \de(\tau)\|_{\dot{H}^s}^2 d\tau\\
 \leq \Bigg[ \|\de(0)\|_{\dot{H}^s}^2 +\frac{C}{\nu_0} \Bigg(\nu_0 \|f^b\|_{L^1 \dot{H}^s} +\|f^l\|_{L^2 \dot{H}^{s-1}}^2 +\|\tUqg\|_{L^\infty\dot{H}^s}^{2(1-\aA)} \|\tUqg\|_{L^2 \dot{H}^{s+1}}^{2\aA} \|\Wet\|_{L^\frac{2}{1-\aA} L^\frac{3}{\aA}}^2\\
 +\|\tUqg\|_{L^\infty \dot{H}^\frac12}^2 \||D|^{s+\aA}\Wet\|_{L^2 L^\frac{6}{1+2\aA}}^2 +\|\tUqg\|_{L^\infty \dot{H}^s}^2 \|\Wet\|_{L^2 \dot{B}_{\frac{3}{\aA},2}^{\aA}}^2 +\||D|^{s+\aB}\Wet\|_{L^2 L^\frac{6}{1+2\aB}}^2\\
 +\|\Wet\|_{L^\infty\dot{H}^s}^{2(1-\aC)} \|\Wet\|_{L^2 \dot{H}^{s+1}}^{2\aC} \|\Wet\|_{L^\frac{2}{1-\aC} L^\frac{3}{\aC}}^2 +\|\Wet\|_{L^\infty \dot{H}^s}^2 \|\Wet\|_{L^2 \dot{B}_{\frac{3}{\aC},2}^{\aC}}^2 \Bigg) \Bigg]\\
\times \exp \frac{C}{\nu_0}\Bigg\{\nu_0 \|f^b\|_{L^1 \dot{H}^s} +\|\n \tUqg\|_{L^2 \dot{H}^{\frac12}}^2 (1+\frac{1}{\nu_0^2}\|\tUqg\|_{L^\infty \dot{H}^{\frac12}}^2) +\frac{1}{\nu_0^{\frac{2\aB}{1-\aB}}}\|\Wet\|_{L^\frac{2}{1-\aB} L^\frac{3}{\aB}}^\frac{2}{1-\aB} +\|\Wet\|_{L^2 \dot{B}_{\frac{3}{\aB},2}^{\aB}}^2 \Bigg\}.
\label{estimCas2}
\end{multline}
Let us precise that compared to \eqref{estimCas1}, the only differences are:
\begin{itemize}
\item the force terms $f^{b,l}$ (dealt as in \cite{FC2, FC3}),
\item the simpler estimates for $F_8$: as precision will be imposed by the truncated terms, we only write:
\begin{multline}
 \Big|(F_8|\de)_{\Hs}\Big| \leq \frac{\nu_0}{16} \|\de\|_{\dot{H}^{s+1}}^2\\
 +\frac{C}{\nu_0} \Big(\|\Wet\|_{L^\infty\dot{H}^s}^{2(1-\aC)} \|\Wet\|_{L^2 \dot{H}^{s+1}}^{2\aC} \|\We\|_{L^\frac{2}{1-\aC} L^\frac{3}{\aC}}^2 +\|\We\|_{L^\infty \dot{H}^s}^2 \|\We\|_{L^2 \dot{B}_{\frac{3}{\aC},2}^{\aC}}^2 \Big).
\end{multline}
\end{itemize}

\subsubsection{Estimates for the truncated quantities}
We will now bound much more precisely than in \cite{FC2, FC3} the external force terms and initial data (see \eqref{f1f2}):
\begin{prop}
 \sl{There exists a constant $\Bo\geq 1$ such that for all $s\in \interv$,
 \begin{equation}
  \begin{cases}
  \vspace{2mm}
   \|f^b\|_{L^1 \Hs } \leq \Bo \left(\ee^{1-2M} +\ee^{M(1-\eta)\delta} +\ee^{\frac{m}6-M(\frac56+\eta \delta)}\right),\\
   \vspace{2mm}
   \|f^l\|_{L^2 \dot{H}^{s-1}} \leq \Bo \left(\ee^{1-2M} +\ee^{M(1-\eta)\delta} +\ee^{m(\frac12 -\eta \delta)}\right),\\
   \|\de(0)\|_{\Hs } \leq \Bo \left( \ee^{\aa_0} +\ee^{1-2M-\gamma} +\ee^{\delta(M-\eta m) -\gamma} +\ee^{\delta\big((\frac12-\eta)m-M\big)-\gamma}\right).
  \end{cases}
\label{estimfext1}
 \end{equation}
 }
\end{prop}
\begin{rem}
 \sl{Note that as we want positive powers of $\ee$, the previous estimates imply the following conditions:
 \begin{equation}
 \begin{cases}
    M,\eta, \eta \delta \in]0,\frac12[,\\
  \eta <\frac{M}{m} <\min(\frac1{5+6\eta \delta}, \frac12-\eta),\\
  \gamma<\min(1-2M, \delta(M-\eta m), \delta \big((\frac12-\eta)m-M\big)).
  \label{CondA}
 \end{cases}
 \end{equation}
 }
\end{rem}

\textbf{Proof :} let us begin with the terms involving $G$: thanks to \eqref{estimvp2}, and Propositions \ref{Gomeganul} and \ref{estimGlb}, we immediately obtain that there exists a constant $\Bo$ (only depending on $\Co, \nu, \nu'$ and $F$) such that for all $s\in \interv$:
$$
\|\cPrr \mathbb{P}_2 G^b\|_{L^1 \Hs } +\|\cPrr \mathbb{P}_2 G^l\|_{L^2 \dot{H}^{s-1}} \leq \Bo \ee \Re^2,\\
$$
Thanks to Lemma \ref{lemaniso} (see appendix), Proposition \ref{estimGlb} and \eqref{estimQGlamb}, the second term in $f_1$ can be bounded (for all $s\in \interv$) according to:
\begin{multline}
\|(Id-\cPrr)G^b\|_{L^1 \Hs} \leq \|(Id-\chi(\frac{|D|}{\Re}))G^b\|_{L^1 \Hs} +\|\chi(\frac{|D|}{\Re}) \chi(\frac{|D_3|}{\re}) G^b\|_{L^1 \Hs}\\
\leq \frac{1}{\Re^{\frac12+\delta-s}}\|(Id-\chi(\frac{|D|}{\Re}))G^b\|_{L^1 \dot{H}^{\frac12 +\delta}} +\Re^s \|\chi(\frac{|D|}{\Re}) \chi(\frac{|D_3|}{\re}) G^b\|_{L^1 L^2}\\
\leq \frac{1}{\Re^{\frac12+\delta-s}}\|G^b\|_{L^1 \dot{H}^{\frac12 +\delta}} +\Re^s (\Re^2 \re)^{\frac23-\frac12}\|\tUqg\cdot \nabla \tUqg\|_{L^1 L^\frac32}\\
\leq \frac{1}{\Re^{\frac12+\delta-s}}\|G^b\|_{L^1 \dot{H}^{\frac12 +\delta}} +\Re^{s+\frac13} \re^{\frac16} \int_0^\infty \|\tUqg(\tau)\|_{L^6} \|\nabla \tUqg (\tau)\|_{L^2} d\tau\\
\leq \Bo (\frac1{\Re^{\frac12+\delta-s}} + \Re^{s+\frac13} \re^{\frac16}).
\end{multline}
which implies the first estimates in \ref{estimfext1} for all $s\in\interv$. Similarly, we have
\begin{equation}
  \|(Id-\chi(\frac{|D|}{\Re}))G^l\|_{L^2 \dot{H}^{s-1}}^2 d\tau \leq \frac{\Bo}{\Re^{2(\frac12+\delta-s)}},
\end{equation}
and using that the expression of $G^l$ (see \eqref{f1f2}) features some derivative $\d_3$, we have for all $s\in \interv$,
\begin{multline}
\|\chi(\frac{|D|}{\Re}) \chi(\frac{|D_3|}{\re}) G^l\|_{L^2 \dot{H}^{s-1}} \leq C_F|\nu-\nu'| \|\chi(\frac{|D|}{\Re}) \chi(\frac{|D_3|}{\re}) \d_3 \nabla \tUqg\|_{L^2 \dot{H}^{s-1}} \\
  \leq C_F|\nu-\nu'| \re^s \|\d_3^{1-s} \nabla \tUqg\|_{L^2 \dot{H}^{s-1}} \leq C_F|\nu-\nu'| \re^s \|\tUqg\|_{L^2 \dot{H}^{1}}.
\end{multline}
Let us now turn to bound the initial data $\de(0)$:
\begin{multline}
 \|\de(0)\|_{\Hs }\leq \|\Uoqg-\tUoqg\|_{\Hs} +\|\cPrr \mathbb{P}_2 \Uoosc\|_{\Hs} + \|(Id-\cPrr) \Uoosc\|_{\Hs}\\
  \leq \Co \ee^{\aa_0} +\|\cPrr \mathbb{P}_2 \Uoosc\|_{\Hs} +\|(Id-\chi(\frac{|D|}{\Re}))\Uoosc\|_{\Hs}+\|\chi(\frac{|D|}{\Re}) \chi(\frac{|D_3|}{\re})\Uoosc\|_{\Hs}.
  \label{estiminit}
\end{multline}
As before, we easily estimate the second and third terms for all $s\in \interv$ by
\begin{equation}
  C_F|\nu-\nu'| \ee \Re^2 \|\Uoosc\|_{\Hs} +\frac{C_F}{\Re^{\frac12 +\delta-s}} \|\Uoosc\|_{\dot{H}^{\frac12 +\delta}} \leq \Bo \ee^{-\gamma} \big[\ee \Re^2 +\frac{1}{\Re^{\frac12+\delta-s}}\big].
\end{equation}
This is here that the $\dot{B}_{q,q}^{\frac12}$-assumption on the initial data will be specifically used (everywhere else we only use the fact that this space is embedded in $\dot{H}^{\frac12-\frac32 \delta}$). To bound the last term, thanks to Proposition \ref{injectionLr} let us write that (we recall that $q=\frac{2}{1+\delta}<2$):
\begin{multline}
 \|\chi(\frac{|D|}{\Re}) \chi(\frac{|D_3|}{\re})\Uoosc\|_{\Hs } =\|\chi(\frac{|D|}{\Re}) \chi(\frac{|D_3|}{\re}) |D|^s \Uoosc\|_{L^2}\\
 \leq C(\Re^2 \re)^{\frac1q-\frac12} \Re^{s-\frac12} \|\chi(\frac{|D|}{\Re}) \chi(\frac{|D_3|}{\re}) |D|^\frac12 \Uoosc\|_{L^q} \leq C \Re^{\delta+s-\frac12} \re^{\frac{\delta}{2}} \||D|^\frac12 \Uoosc\|_{\dot{B}_{q,q}^0}\\
 \leq C\Re^{\delta+s-\frac12} \re^{\frac{\delta}{2}} \|\Uoosc\|_{\dot{B}_{q,q}^{\frac12}}.
\end{multline}
note that this can be done only if $s\geq \frac12$. In the case $s\in [\frac12-\eta \delta, \frac12[$, we simply go back to \eqref{estiminit} and write that (taking advantage of the frequency localization):
\begin{multline}
 \|(Id-\cPrr)\Uoosc\|_{\Hs} \leq \frac1{\re^{\frac12 -s}} \|(Id-\cPrr)\Uoosc\|_{\dot{H}^\frac12}\\
 \leq \frac{C_F}{\re^{\frac12 -s}} \left(\frac1{\Re^{\delta}} \|\Uoosc\|_{\dot{H}^{\frac12 +\delta}} +\Re^{\delta} \re^{\frac{\delta}{2}} \|\Uoosc\|_{\dot{B}_{q,q}^{\frac12}} \right).
\end{multline}
We can sum up as follows: for all $s\in \interv$
\begin{multline}
 \|(Id-\cPrr)\Uoosc\|_{\Hs} \leq \Co \ee^{-\gamma}
 \times \begin{cases}
         \frac1{\Re^{(1-\eta)\delta}} +\Re^{(1+\eta)\delta} \re^{\frac{\delta}2} & \mbox{ if } s\in[\frac12, \frac12+\eta \delta],\\
         \frac1{\re^{\eta \delta}} (\frac1{\Re^\delta}+\Re^\delta \re^{\frac{\delta}2})& \mbox{ if } s\in[\frac12 -\eta\delta, \frac12]
        \end{cases}\\
        \leq \Co \ee^{-\gamma} \times \begin{cases}
         \ee^{M \delta(1-\eta)} +\ee^{\delta(\frac{m}2 -(1+\eta)M)} & \mbox{ if } s\in[\frac12, \frac12+\eta \delta],\\
         \ee^{\delta(M-m\eta)} +\ee^{\delta\big((\frac12-\eta)m-M\big)} & \mbox{ if } s\in[\frac12 -\eta\delta, \frac12].
        \end{cases}
\end{multline}
As
$$
\begin{cases}
 M (1-\eta)-(M-m\eta)=\eta(m-M),\\
 (\frac{m}2 -(1+\eta)M) -\big((\frac12-\eta)m-M\big)= \eta(m-M),
\end{cases}
$$
and as $m>M$ (see \eqref{CondA}), we obtain the announced result. $\blacksquare$

\subsubsection{Strichartz estimates for $\Wet$}

We will need the following Strichartz estimates to complete our bootstrap argument:
\begin{prop}
 \sl{There exist $\ee_0, \Bo>0$ such that for any $\aa>0$ and $\ee<\ee_0$, $\Wet$ satisfies:
 \begin{equation}
  \begin{cases}
   \|\Wet\|_{\tilde{L}^2 \dot{B}_{\frac3{\aa}, 2}^\aa} \leq \Bo \ee^{\frac14 -\frac{\aa}3 -M(\frac92-4\aa-\delta) -m(\frac92-2\aa)} \leq \Bo \ee^{\frac14 -\frac{\aa}3 -\frac92(M+m)},\\
   \|\Wet\|_{L^\frac2{1-\aa} L^\frac3{\aa}} \leq \Bo \ee^{\frac14 -\frac{\aa}3 -M(\frac92-3\aa-\delta) -m(\frac92-3\aa)} \leq \Bo \ee^{\frac14 -\frac{\aa}3 -\frac92(M+m)},\\
   \||D|^{s+\aa}\Wet\|_{L^2 L^\frac6{1+2\aa}}\leq \Bo \ee^{\frac1{12} -\frac{\aa}3 -M(\frac72-3\aa) -m(\frac72-2\aa)} \leq \Bo \ee^{\frac1{12} -\frac{\aa}3 -\frac72(M+m)}.
  \end{cases}
\label{StriBB}
 \end{equation}
 }
\end{prop}
\textbf{Proof:} using Proposition \ref{Estimdispnu2} in the case $(d,p,r,q)=(\aa, 2, \frac3{\aa}, 2)$, we obtain that
\begin{multline}
 \|\Wet\|_{\tilde{L}^2 \dot{B}_{\frac3{\aa}, 2}^\aa} \leq \Bo \ee^{\frac14(1-\frac{4\aa}3)} \frac{\Re^{4-3\aa}}{\re^{\frac72-2\aa}}\\
 \times \left( \|\cPrr \Uoosc\|_{\dot{H}^\aa}+ \|\cPrr G^b\|_{L^1 \dot{H}^\aa} +\frac1{\nu_0^\frac12 \re} \|\cPrr G^l\|_{L^2 \dot{H}^\aa}\right)\\
 \leq \Bo \ee^{\frac14(1-\frac{4\aa}3)} \frac{\Re^{4-3\aa}}{\re^{\frac72-2\aa}} \times \left( \frac1{\re^{\frac12-\frac{3\delta}2-\aa}}\|\Uoosc\|_{\dot{H}^{\frac12 -\frac{3\delta}2}}+ \|G^b\|_{L^1 \dot{H}^\aa} +\frac1{\nu_0^\frac12 \re} \Re^{\frac12-\delta-\aa}\|G^l\|_{L^2 \dot{H}^{\frac12+\delta}}\right)\\
 \leq \Bo \ee^{\frac14-\frac{\aa}3-M(4-3\aa)-m(\frac72-2\aa)}\left(\ee^{-\gamma-m(\frac12-\frac{3\delta}2-\aa)} +\ee^{-m-M(\frac12-\delta-\aa)}\right).
\end{multline}
From \eqref{CondA}, we know that $\gamma<\delta M$ so that
$$
m+M(\frac12-\delta-\aa)-\left( \gamma+m(\frac12-\frac{3\delta}2-\aa)\right)=M(\frac12 -\delta -\aa)+m(\frac12 +\frac{3\delta}{2}+\aa) -\gamma>0,
$$
which leads to the first estimate. Similarly, considering Proposition \ref{Estimdispnu2} in the case $(d,p,r,q)=(0, \frac2{1-\aa}, \frac3{\aa}, 2)$, we get 
(thanks to proposition \ref{injectionLr}):
\begin{multline}
 \|\Wet\|_{L^\frac2{1-\aa} L^{\frac3{\aa}}} \leq \|\Wet\|_{\tilde{L}^2 \dot{B}_{\frac3{\aa}, 2}^0} \leq \Bo \ee^{\frac14(1-\frac{4\aa}3)} \frac{\Re^{4-3\aa}}{\re^{\frac72-3\aa}}\\
 \times \left( \frac1{\re^{\frac12-\frac{3\delta}2}}\|\Uoosc\|_{\dot{H}^{\frac12 -\frac{3\delta}2}}+ \|G^b\|_{L^1 L^2} +\frac1{\nu_0^\frac12 \re} \Re^{\frac12-\delta}\|G^l\|_{L^2 \dot{H}^{\frac12+\delta}}\right)\\
 \leq \Bo \ee^{\frac14-\frac{\aa}3-M(4-3\aa)-m(\frac72-3\aa)}\left(\ee^{-\gamma-m(\frac12-\frac{3\delta}2)} +\ee^{-m-M(\frac12-\delta)}\right),
\end{multline}
which leads to the second estimates. In the case $(d,p,r,q)=(s+\aa, 2, \frac6{1+2\aa}, 2)$, we obtain that (provided that $0<\aa<\delta+\frac12-s$):
\begin{multline}
 \||D|^{s+\aa}\Wet\|_{L^2 L^{\frac6{1+2\aa}}} \leq \|\Wet\|_{\tilde{L}^2 \dot{B}_{\frac6{1+2\aa}, 2}^{s+\aa}} \leq \Bo \ee^{\frac1{12}(1-4\aa)} \frac{\Re^{\frac52-3\aa}}{\re^{\frac52-2\aa}}\\
 \times \left( \|\cPrr \Uoosc\|_{\dot{H}^{s+\aa}}+ \|\cPrr G^b\|_{L^1 \dot{H}^{s+\aa}} +\frac1{\nu_0^\frac12 \re} \|\cPrr G^l\|_{L^2 \dot{H}^{s+\aa}}\right)\\
 \leq \Bo \ee^{\frac1{12}(1-4\aa)} \frac{\Re^{\frac52-3\aa}}{\re^{\frac52-2\aa}} \times \left( \ee^{-\gamma} +1 +\frac1{\re} \Re^{s+\aa+\frac12-\delta}\|\cPrr \tUqg\|_{L^2 \dot{H}^{\frac32+\delta}}\right)\\
 \leq \Bo \ee^{\frac1{12}-\frac{\aa}3-M(\frac52-3\aa)-m(\frac52-2\aa)}\left(\ee^{-\gamma} +\ee^{-m-M}\right),
\end{multline}
which concludes the proof. $\blacksquare$

\subsubsection{Bootstrap}

We are now able to conclude the bootstrap argument (see previous section and \cite{FC2, FC3}). Gathering \eqref{estimCas2}, \eqref{estimfext1} and \eqref{StriBB}, we obtain that for all $t\leq T_\ee$,
\begin{multline}
 \|\de(t)\|_{\dot{H}^s}^2 +\frac{\nu_0}{2} \int_0^t \|\n \de(\tau)\|_{\dot{H}^s}^2 d\tau \leq \Do \Bigg[\ee^{2\aa_0} +\ee^{2(1-2M-\gamma)} +\ee^{2\big(\delta(M-\eta m) -\gamma\big)}\\
 +\ee^{2\Big(\delta\big((\frac12-\eta)m-M\big)-\gamma\Big)} +\ee^{1-2M} +\ee^{M(1-\eta)\delta} +\ee^{\frac{m}6-M(\frac56+\eta \delta)} +\ee^{2(1-2M)} +\ee^{2M(1-\eta)\delta}\\
 +\ee^{2m(\frac12 -\eta \delta)} +\ee^{\frac14-\frac{\aa}3-\frac92(M+m)-\gamma} +\ee^{\frac1{12} -\frac{\aa}3 -\frac72(M+m)} \Bigg] \times \exp \frac{C}{\nu_0}\Bigg\{1 +\ee^{1-2M} +\ee^{M(1-\eta)\delta}\\
 +\ee^{\frac{m}6-M(\frac56+\eta \delta)} +\ee^{\min(2,\frac2{1-\aa})\left(\frac14 -\frac{\aa}3 -\frac92(M+m)\right)} \Bigg\}.
\label{estimCas2b}
\end{multline}
For more simplicity we will ask, instead of the second condition from \eqref{CondA}, that:
$$
2\eta \leq \frac{M}{m} \leq \frac12 \min(\frac1{5+6\eta \delta}, \frac12-\eta).
$$
This obviously implies that $\eta\leq \frac1{10}$, so we will finally ask that:
\begin{equation}
 \begin{cases}
    M\in ]0,\frac14],\mbox{ }\eta \in]0,\frac1{10}[,\\
  2\eta \leq\frac{M}{m} \leq \frac12 \frac1{5+\delta},\\
  \gamma<\min(\frac12(1-2M), \frac12 \delta(M-\eta m)), \frac12\delta \big((\frac12-\eta)m-M\big)).
  \label{CondB}
 \end{cases}
 \end{equation}
Moreover, if we take $\aa=\gamma$ and ask that:
\begin{equation}
 \begin{cases}
  \frac92 (M+m) \leq \frac18 & \mbox{ and }\frac43 \delta \leq \frac12(\frac14-\frac92 (M+m)),\\
  \frac72 (M+m) \leq \frac1{24} & \mbox{ and }\frac{\delta}3 \leq \frac12(\frac1{12}-\frac72 (M+m)).
 \end{cases}
\end{equation}
 As $M\leq \frac{m}{10}$, this is realized when
$$
 \begin{cases}
    m\in ]0,\frac1{100}],\\
  2\eta \leq\frac{M}{m} \leq \frac12 \frac1{5+\delta}.
 \end{cases}
$$
When
$$
\gamma \leq \min(\frac{M\delta}{4}, \frac{m\delta}{16}, \frac{m}{12}, \frac1{32}) =\frac{M\delta}{4},
$$
we obtain that all power of $\ee$ in the exponential is positive son that for small enough $\ee$, we get that for all $s\in \interv$ and $t\leq T_\ee$:
\begin{equation}
 \|\de(t)\|_{\dot{H}^s}^2 +\frac{\nu_0}{2} \int_0^t \|\n \de(\tau)\|_{\dot{H}^s}^2 d\tau \leq \Do e^{2\Do} \ee^{\min(2\aa_0, \frac{M\delta}2)},
\label{estimfinaleB}
 \end{equation}
so that we finally end-up with (for small enough $\ee$), $\de(T_\ee) \leq \frac{\nu_0}{8C}$ which clearly contradicts the maximality of $T_\ee$. We can conclude that $T_\ee=T_\ee^*$ and then the previous estimate is valid for all $t<T_\ee^*$, which implies for $s=\frac12$ that the integral in \eqref{critereexpl} is finite. Therefore $T_\ee^*=\infty$ and \eqref{estimfinaleB} is then valid for all $t\geq 0$. The rest of the theorem is done as for the case $\nu=\nu'$. $\blacksquare$

\section{Appendix}

\subsection{Notations and Sobolev spaces}

For $s\in\R$, $\Hs $ and $H^s$ are the classical homogeneous/inhomogeneous Sobolev spaces in $\R^3$ endowed with the norms:
$$
\|u\|_{\Hs }^2=\int_{\R^3} |\xi|^{2s} |\hat{u}(\xi)| d\xi, \quad \mbox{and} \quad \|u\|_{H^s}^2=\int_{\R^3} (1+|\xi|^2)^s |\hat{u}(\xi)| d\xi.
$$
We also use the following notations: if $E$ is a Banach space and $T>0$,
$$
\cC_T{E} =\cC([0,T], E), \quad \mbox{and} \quad L_T^p{E} =L^p([0,T], E).
$$
Let us recall the Sobolev injections, and product laws:
\begin{prop}
 \sl{There exists a constant $C>0$ such that if $s<\frac{3}{2}$, then for any $u\in\Hs (\R^3)$, $u\in L^{p}(\R^3)$ with $p=\frac{6}{3-2s}$ and
 $$
 \|u\|_{L^{p}}\leq C \|u\|_{\Hs }.
 $$}
\end{prop}
\begin{prop} (\cite{Dbook}, chapter 2)
 \sl{There exists a constant $C$ such that for any $(u,v)\in \dot{H}^{s_1}(\R^3) \times \dot{H}^{s_2}(\R^3)$, if $s_1,s_2\in ]-\frac{3}{2},\frac{3}{2}[$ and $s_1+s_2>0$ then $uv \in \dot{H}^{s_1+s_2-\frac{3}{2}}(\R^3)$ and we have:
 $$
 \|uv\|_{\dot{H}^{s_1+s_2-\frac{3}{2}}} \leq C \|u\|_{\dot{H}^{s_1}} \|v\|_{\dot{H}^{s_2}}.
 $$
 }
 \label{HProd}
\end{prop}

\subsection{Besov spaces}

We refer to Chapter 2 from \cite{Dbook} for an in-depth presentation of the classical homogeneous and inhomogeneous Besov and Sobolev spaces. We also refer to the appendix of \cite{FCpochesLp} for a quick presentation.

Let us just recall that $\psi$ a smooth radial function supported in the ball $B(0,\frac{4}{3})$, equal to 1 in a neighborhood of $B(0,\frac{3}{4})$ and such that $r\mapsto \psi(r.e_r)$ is nonincreasing over $\R_+$. If we set $\varphi(\xi)=\psi(\xi/2)-\psi(\xi)$, then $\varphi$ is compactly supported in the annulus $\cC=\{\xi\in \R^d, c_0=\frac{3}{4}\leq |\xi|\leq C_0=\frac{8}{3}\}$ and we define the homogeneous dyadic blocks: for all $j\in \Z$,
$$
\ddj u:= \varphi(2^{-j}D)u =2^{jd} h(2^j.)* u, \quad \mbox{with } h=\cF^{-1} \varphi.
$$
We recall that $\hat{\phi(D)u}(\xi)=\phi(\xi) \hat{u} (\xi)$ and we can define the homogeneous Besov norms and spaces:
\begin{defi}
\sl{For $s\in\R$ and $1\leq p,r\leq\infty,$ we set
$$
\|u\|_{\dot B^s_{p,r}}:=\bigg(\sum_{l\in \Z} 2^{rls}
\|\ddl u\|^r_{L^p}\bigg)^{\frac{1}{r}}\ \text{ if }\ r<\infty
\quad\text{and}\quad
\|u\|_{\dot B^s_{p,\infty}}:=\sup_{l} 2^{ls}
\|\ddl u\|_{L^p}.
$$
The homogeneous Besov space $\dot B^s_{p,r}$ is the subset of tempered distributions such that $\lim_{j \rightarrow -\infty} \|\dot{S}_j u\|_{L^\infty}=0$ and $\|u\|_{\dot B^s_{p,r}}$ is finite (where $\dot{S}_j u=\Sum_{l\leq j-1} \ddl u=\psi(2^{-j}D)u$).
}
\end{defi}
\begin{itemize}
\item The space $\dot B^s_{p,r}$ is complete whenever
$s<d/p,$ or $s\leq d/p$ and $r=1$.
  \item For any $p\in[1,\infty],$ we have the  continuous embedding $\dot B^0_{p,1}\hookrightarrow L^p\hookrightarrow \dot B^0_{p,\infty}.$
\item If  $\sigma\in\R,$ $1\leq p_1\leq p_2\leq\infty$ and $1\leq r_1\leq r_2\leq\infty,$
  then $\dot B^{\sigma}_{p_1,r_1}\hookrightarrow
  \dot B^{\sigma-d(\frac1{p_1}-\frac1{p_2})}_{p_2,r_2}.$
  \item The space  $\dot B^{\frac dp}_{p,1}$ is continuously embedded in   the set  of
bounded  continuous functions (going to $0$ at infinity if    $p<\infty$).
\item $\dot{H}^s=\dot{B}_{2,2}^s$.
\item Interpolation: if $1\leq p,r_1,r_2,r\leq\infty,$ $\sigma_1\not=\sigma_2$ and $\theta\in(0,1)$:
  \begin{equation}\label{R-E9}
  \|f\|_{\dot B^{\theta\sigma_2+(1-\theta)\sigma_1}_{p,r}}\lesssim\|f\|_{\dot B^{\sigma_1}_{p,r_1}}^{1-\theta}
  \|f\|_{\dot B^{\sigma_2}_{p,r_2}}^\theta.
  \end{equation}
  \end{itemize}
\begin{prop}
 \sl{\cite{Dbook} We have the following continuous injections:
$$
 \begin{cases}
\mbox{For any } p\geq 1, & \dot{B}_{p,1}^0 \hookrightarrow L^p,\\
\mbox{For any } p\in[2,\infty[, & \dot{B}_{p,2}^0 \hookrightarrow L^p,\\
\mbox{For any } p\in[1,2], & \dot{B}_{p,p}^0 \hookrightarrow L^p.
\end{cases}
$$
}
 \label{injectionLr}
\end{prop}
Let us then define the spaces $\tilde L^\rho_T \dot B^s_{p,r}$ from the following norm:
\begin{defi}\label{d:espacestilde}
For $T>0,$ $s\in\R$ and  $1\leq r,\rho\leq\infty,$
 we set
$$
\|u\|_{\tilde L_T^\rho \dot B^s_{p,r}}:=
\bigl\Vert2^{js}\|\ddq u\|_{L_T^\rho L^p}\bigr\Vert_{\ell^r(\Z)}.
$$
\end{defi}
Any product of two distributions $u$ and $v$ may be formally written through the Bony decomposition:
\begin{equation}\label{eq:bony}
uv=T_uv+T_vu+R(u,v), \mbox{ where}
\end{equation}
$$
T_uv:=\sum_l\dot S_{l-1}u\ddl v,\quad
T_vu:=\sum_l \dot S_{l-1}v\ddl u\ \hbox{ and }\ 
R(u,v):=\sum_l\sum_{|l'-l|\leq1}\ddl u\,\dot\Delta_{l'}v.
$$
The above operator $T$ is called a ``paraproduct'' whereas $R$ is called a ``remainder''. We refer to \cite{Dbook} for general properties and for paraproduct and remainder estimates but we can recall that (if $\frac1{r}=\frac1{r_1}+\frac1{r_2}$ and $\frac1{p}=\frac1{p_1}+\frac1{p_2}$):
\begin{itemize}
 \item For any $s\in \R$, $\|T_u v\|_{\dot{B}_{p,r}^s} \lesssim \|u\|_{L^\infty} \|v\|_{\dot{B}_{p,r}^s}$,
 \item For any $(s,t)\in \R_-^*\times \R$, $\|T_u v\|_{\dot{B}_{p,r}^{s+t}} \lesssim \|u\|_{\dot{B}_{p_1,r_1}^s} \|v\|_{\dot{B}_{p_2,r_2}^t}$, 
 \item For any $s,t\in \R$ with $s+t>0$, $\|R(u,v)\|_{\dot{B}_{p,r}^{s+t}} \lesssim \|u\|_{\dot{B}_{p_1,r_1}^s} \|v\|_{\dot{B}_{p_2,r_2}^t}$.
\end{itemize}

\subsection{Dispersion and Strichartz estimates}

Consider the following system:
\begin{equation}
\begin{cases}
\d_t f-(L-\frac{1}{\ee} \mathbb{P} \cA) f=\Fe,\\
f_{|t=0}=f_0.
\end{cases}
\label{systdisp}
\end{equation}
If we apply the Fourier transform, the equation becomes (see \cite{FC} for precisions):
$$
\d_t \hat{f}- \mathbb{B}(\xi, \ee)\hat{f}=\hat{\Fe},
$$
where
$$\mathbb{B}(\xi, \ee)= \hat{L-\frac{1}{\ee} \mathbb{P} \cA} =\left(
\begin{array}{cccc}
\displaystyle{-\nu|\xi|^2+\frac{\xi_1\xi_2}{\ee
  |\xi|^2}} & \displaystyle{\frac{\xi_2^2+\xi_3^2}{\ee
  |\xi|^2}} & 0 & \displaystyle{\frac{\xi_1\xi_3}{\ee F |\xi|^2}}\\
\displaystyle{-\frac{\xi_1^2+\xi_3^2}{\ee
  |\xi|^2}} & \displaystyle{-\nu|\xi|^2-\frac{\xi_1\xi_2}{\ee
  |\xi|^2}} & 0 & \displaystyle{\frac{\xi_2\xi_3}{\ee F |\xi|^2}}\\
\displaystyle{\frac{\xi_2\xi_3}{\ee |\xi|^2}} &
\displaystyle{-\frac{\xi_1\xi_3}{\ee
  |\xi|^2}} & \displaystyle{-\nu |\xi|^2} & \displaystyle{-\frac{\xi_1^2+\xi_2^2}{\ee F
  |\xi|^2}}\\
0 & 0 & \displaystyle{\frac{1}{\ee F}} &
\displaystyle{-\nu'|\xi|^2}
\end{array}
\right).
$$
For $0<r<R$ we will denote by $\cC_{r,R}$ the following set:
$$
\cC_{r,R} =\{\xi \in \R^3, \quad |\xi|\leq R \mbox{ and } |\xi_3|\geq r\}.
$$
We also introduce the following frequency truncation operator on $\cC_{r,R}$:
\begin{equation}
 \cP_{r,R}=\chi (\frac{|D|}{R})\big(1-\chi (\frac{|D_3|}{r})\big),
\label{PrR}
\end{equation}
where $\chi$ is the smooth cut-off function introduced before and ($\mathcal{F}^{-1}$ is the inverse Fourier transform):
$$
\chi (\frac{|D|}{R}) f=\mathcal{F}^{-1} \Big(\chi(\frac{|\xi|}{R}) \hat{f}(\xi)\Big) \quad \mbox{and} \quad \chi (\frac{|D_3|}{r}) f=\mathcal{F}^{-1} \Big(\chi(\frac{|\xi_3|}{r}) \hat{f}(\xi)\Big),
$$
and $|D|^s$ the classical derivation operator: $|D|^s f =\mathcal{F}^{-1} (|\xi|^s \hat{f}(\xi)).$

In what follows we will use it for particular radii $r_\ee=\ee^m$ and $R_\ee =\ee^{-M}$, where $m$ and $M$ will be precised later. Let us end with the following anisotropic Bernstein-type result (we refer to \cite{FC}, and to \cite{Dragos1} for more general anisotropic estimates):
\begin{lem}
\sl{There exists a constant $C>0$ such that for all function $f$, $\aa>0$, $1\leq q \leq p \leq \infty$ and all $0<r<R$, we have
\begin{equation}
\begin{cases}
\vspace{1mm}
\displaystyle{\|\chi (\frac{|D|}{R}) \chi (\frac{|D_3|}{r}) f\|_{L^p} \leq C \|f\|_{L^p},}\\
\displaystyle{\|\chi (\frac{|D|}{R}) \chi (\frac{|D_3|}{r}) f\|_{L^p} \leq C(R^2 r)^{\frac{1}{q}-\frac{1}{p}} \|\chi (\frac{|D|}{R}) \chi (\frac{|D_3|}{r}) f\|_{L^q}.}
\end{cases}
\end{equation}
Moreover if $f$ has its frequencies located in $\cC_{r,R}$, then
$$
\||D|^\aa f\|_{L^p} \leq C R^\aa \|f\|_{L^p}. \blacksquare
$$
}
\label{lemaniso}
\end{lem}

\subsubsection{Eigenvalues, projectors}

We begin with the eigenvalues and eigenvectors of matrix $\mathbb{B}(\xi, \ee)$. We refer to \cite{FC, FC2, FC3, FC5, FCpochesLp} for details about the following proposition. We will only state the results and skip details as the proof is an adaptation of Proposition $3.1$ from \cite{FC5} (there in the anisotropic case).

\begin{prop}
\label{estimvp}
\sl{
If $\nu\neq \nu'$ there exists $\ee_0>0$ such that for all $\ee< \ee_0$, for all $\re=\ee^m$ and $\Re =\ee^{-M}$, with $M<1/4$ and $3M+m<1$, and for all $\xi \in \mathcal{C}_{r_{\ee}, R_{\ee}}$, the matrix $\mathbb{B}(\xi, \ee) = \widehat{L-\frac{1}{\ee} \mathbb{P} \mathcal{A}}$ is diagonalizable and its eigenvalues have the following asymptotic expansions with respect to $\ee$:
\begin{equation}
\label{vp}
\begin{cases}
\vspace{0.2cm} \mu_0 = -\nu |\xi|^2,\\
\vspace{0.2cm} \mu = -(\nu\xi_1^2+ \nu\xi_2^2 + \nu'F^2\xi_3^2)\frac{|\xi|^2}{|\xi|_F^2} + \ee^2 D(\xi, \ee),\\
\vspace{0.2cm} \lambda = -\tau (\xi)|\xi|^2+i\frac{|\xi|_F}{\ee F   |\xi|}+ \ee E(\xi, \ee),\\
\vspace{0.2cm} \overline{\lambda} = -\tau (\xi)|\xi|^2-i\frac{|\xi|_F}{\ee F |\xi|}+ \ee \bar{E}(\xi, \ee),
\end{cases}
\end{equation}
where $|\xi|_F^2 = \xi_1^2 + \xi_2^2 + F^2\xi_3^3$, and $D,E$ denote remainder terms satisfying for all $\xi \in\cC_{\re, \Re}$:
$$
\begin{cases}
\vspace{1mm}
\ee^2 |D(\xi, \ee)|\leq C_F |\nu-\nu'|^3 \ee^2 |\xi|^6\leq C_F |\nu-\nu'|^3\ee^{2-6M} \ll 1,\\
\ee |E(\xi, \ee)|\leq C_F |\nu-\nu'|^2 \ee |\xi|^4 \leq C_F |\nu-\nu'|^2\ee^{1-4M} \ll 1,\\
\ee |\partial_{\xi_2} E(\xi, \ee)|\leq C_F |\nu-\nu'|^2 \ee |\xi|^3 \leq C_F |\nu-\nu'|^2\ee^{1-3M} \ll 1,\\
\end{cases}
$$
and
$$
\tau(\xi)=\frac{\nu}{2}\Big(1+\frac{F^2 \xi_3^2}{|\xi|_F^2}\Big) +\frac{\nu'}{2}\Big(1-\frac{F^2   \xi_3^2}{|\xi|_F^2}\Big) \geq \min (\nu,\nu')>0.
$$
Moreover, if we denote by $\mathcal{P}_i(\xi, \ee)$, the projectors onto the eigenspaces corresponding to $\mu$, $\lambda$ and $\overline{\lambda}$ ($i\in\{2,3,4\}$), and set
\begin{equation}
\label{projsevp} \mathbb{P}_i(u)=\mathcal{F}^{-1}\bigg(\mathcal{P}_i\big(\xi, \ee)(\widehat{u}(\xi)\big)\bigg),
\end{equation}
then for any divergence-free vector field $f$ whose Fourier transform is supported in $\mathcal{C}_{\re, \Re}$ and $s\in\R$, we have the following estimates:
\begin{equation}
\label{estimvp2}
\|\mathbb{P}_2 f\|_{\Hs } \leq C_F \|f\|_{\Hs } \times
\begin{cases}
{1} & {\mbox{if} \quad \Omega(f)\neq 0,}\\
{|\nu-\nu'| \ee \Re^2 =|\nu-\nu'| \ee^{1-2M}} & {\mbox{if} \quad \Omega(f)= 0,}
\end{cases}
\end{equation}
and for $i=3,4$,
\begin{equation}
\label{estimvp34}
\|\mathbb{P}_i f\|_{\Hs }\leq C_F \frac{\Re}{\re} \|f\|_{\Hs } =C_F \ee^{-(m+M)} \|f\|_{\Hs }.
\end{equation}
Finally, if we define $\mathbb{P}_{3+4} f \overset{def}{=} (\mathbb{P}_3 +\mathbb{P}_4) f= (I_d-\mathbb{P}_2) f$ (as $\div f=0$), then
\begin{equation}
\label{estimvp3+4}
\|\mathbb{P}_{3+4} f\|_{\Hs }\leq C_F (1+ |\nu-\nu'| \ee \Re^2) \|f\|_{\Hs }.
\end{equation}
}
\end{prop}
\begin{rem}
 \sl{In the case $\nu=\nu'$ everything is simpler: the eigenvalues have simple explicit expressions: $-\nu |\xi|^2$ (double, $\mu$ and $\mu_0$ coincide), $-\nu|\xi|^2\pm \frac{i}{\ee}\frac{|\xi|_F}{F|\xi|}$, the eigenvectors do not depend on $\ee$ and are mutually orthogonal (so that $\mathbb{P}_i$ are of norm $1$) and this basis exactly correspond to the $QG/osc$ decomposition (for divergence-free vector fields): $\cP=\mathbb{P}_{3+4}$ and $\cQ=\mathbb{P}_2$ so that the quasigeostrophic part only depends on $W_2$ and the oscillating part only depends on $W_{3,4}$. Finally the operator $\G$ reduces to $\nu \Delta$. We refer to \cite{FC3} (Appendix B) or \cite{FCestimLp, FCpochesLp} for more details.
 }
 \label{remnunu}
\end{rem}

\begin{rem}
\sl{We emphasize that the higher order term $\mu$ is the Fourier symbol of the quasi-geostrophic operator $\G$. Moreover, the dispersion is related to the term $i\frac{|\xi|_F}{\ee F |\xi|}$, and when $F=1$ this term reduces to the constant $\frac{i}{\ee}$. This is why dispersion does not occur in the case $F=1$ (we refer to \cite{Chemin2, FCF1} for a study of the asymptotics in the special cas $F=1$).}
\end{rem}

\subsubsection{Dispersion, Strichartz estimates}

Combining Proposition 3 from \cite{FCpochesLp} (covering the range $p\geq 4$) with the convolution arguments from the appendix of \cite{FC2} allows us to cover the full range $p\geq 1$ and obtain the following Strichartz estimates satisfied by the last two projections of the solution of System\eqref{systdisp}:
\begin{prop}
\sl{Assume that $f$ satisfies \eqref{systdisp} on $[0,T[$ where $\div f_0=0$ and the frequencies of $f_0$ and $F$ are localized in $\cC_{\re, \Re}$. Then there exists a constant $C=C_{F,p, \nu, \nu'}>0$ such that for $i\in\{3,4\}$ and $p\geq 1$, we have
$$
\|\mathbb{P}_i f\|_{L_T^p L^\infty} \leq C K(\ee) \left(\|f_0\|_{L^2} +\int_0^T \|\Fe(\tau)\|_{L^2} d\tau\right).
$$
where
$$
K(\ee)=
\begin{cases}
\vspace{2.5mm}
 \displaystyle{\ee^{\frac 14} \frac{\Re^4}{\re^{\frac 52 + \frac 2p}} \big[\frac{4}{\nu_0}(\frac1p -\frac14)\big]^{\frac1p -\frac14}  =\ee^{\frac 14-\left(4M +(\frac 52+\frac 2p)m\right)} \big[\frac{4}{\nu_0}(\frac1p -\frac14)\big]^{\frac1p -\frac14},} & {\mbox{if }p\in [1,4],}\\
 \displaystyle{\ee^{\frac 1p} \frac{\Re^{\frac 52+ \frac 6p}}{\re^{2 + \frac 4p}} =\ee^{\frac 1p-\left((\frac{5}{2}+\frac{6}{p})M +(2+\frac{4}{p})m\right)},} & {\mbox{if }p\geq 4.}
\end{cases}
$$
\label{Strichartz1}
}
\end{prop}
Unfortunately these estimates would be completely useless in our case: we need more flexibility than only $L^p-L^\infty$-estimates, and in the case $\nu\neq \nu'$ we need to take into account the second term $G^l$ as done in \cite{FC2}. We begin with the case $\nu=\nu'$ where we have to deal with the fact that we obtain Strichartz estimates on $\We$ which is not frequency localized (we improve the method from \cite{FC3} Appendix B). Then we deal with the case $\nu\neq \nu'$.

\subsubsection{Strichartz estimates in the case $\nu=\nu'$}

The main result of this section is stated as follows:

\begin{prop}
 \sl{There exists a constant $C_F>0$, such that for any $d\in \R$, $r>4$, $q\geq 1$ and
 $$
 \theta\in ]0,\frac{\frac12 -\frac1r}{1-\frac{4}{r}}[\cap]0,1], \quad p\in[1, \frac{4}{\theta (1-\frac{4}{r})}],
 $$ 
 if $f$ solves \eqref{systdisp} for initial data $f_0$ and external force $\Fe$ both with zero divergence and potential vorticity, then ($c_0$ refers to the smaller constant appearing in the Littlewood-Paley decomposition, usually $c_0=\frac34$.)
 \begin{equation}
  \||D|^d f\|_{\tilde{L}_t^p\dot{B}_{r, q}^0} \leq C_F\frac{C_{p,\theta,r}}{\nu^{\frac{1}{p}-\frac{\theta}{4}(1-\frac{4}{r})}} \ee^{\frac{\theta}{4}(1-\frac{4}{r})} \times \left( \|f_0\|_{\dot{B}_{2, q}^\sigma} +\int_0^t  \|\Fe(\tau)\|_{\dot{B}_{2, q}^\sigma} d\tau \right),
 \end{equation}
 where
 $$
 \begin{cases}
  \sigma= d+\frac32-\frac{3}{r}-\frac{2}{p}+\frac{\theta}{2} (1-\frac{4}{r}),\\
  C_{p,\theta,r}=\left[\frac{2}{c_0^2}\Big(\frac{1}{p}-\frac{\theta}{4}(1-\frac{4}{r})\Big)\right]^{\frac{1}{p}-\frac{\theta}{4}(1-\frac{4}{r})} \displaystyle{\frac{2^{\frac12 \big(1-\frac{2}{r}-2\theta (1-\frac{4}{r})\big)}}{1-2^{-\frac12 \Big(1-\frac{2}{r}-2\theta (1-\frac{4}{r})\Big)}}}.
 \end{cases}
 $$
\label{Estimdispnu}
 }
\end{prop}
\begin{rem}
\sl{It is interesting to compare our Strichartz estimates with the ones from \cite{IMT} (see Proposition \ref{StriIMT}). In our estimates we use the range $r>4$ whereas in Proposition \ref{StriIMT} is considered the case $r\in]2,4[$ and they use it for $r$ close to $3$. Our index $p$ is mostly equal to $2$ but we can reach $p=1$ (which is useful when there are derivatives), whereas in \cite{IMT}, $p> \frac1{1-\frac2{r}}>2$. 
.}
\label{RqStri}
\end{rem}

\textbf{Proof: } Let us first assume that $\Fe=0$. As $\nu=\nu'$, the fact that $f_0$ is divergence-free and with zero potential vorticity implies that:
$$
f_0=\mathbb{P} f_0=\cP \mathbb{P} f_0= \mathbb{P}_{3+4} \mathbb{P} f_0 =\mathbb{P}_{3+4} f_0,
$$
So that we only consider the last two eigenvalues (we recall the eigenvectors are orthogonal). The idea is here to push further the Strichartz estimates without frequency truncation we obtained in \cite{FC3}: we will once more use a simple non-stationnary phase argument (see for example the works of Chemin, Desjardins, Gallagher and Grenier, we refer to \cite{CDGG, CDGG2, CDGGbook}). As outlined previously, in this special case there is no need to truncate in frequency through the operator $\cPrr$ but within the computations we will truncate considering the vertical Littlewood-Paley decomposition ($\ddk u=\varphi(2^{-j} D_3)u$):
$$
\|\ddj f\|_{L_t^p L_x^r}= \|\ddj f\|_{L^p L^r} \leq \Sum_{k=-\infty}^{j+1} \|\ddk \ddj f\|_{L^p L^r}.
$$
Now we will use the methods leading to the general Strichartz estimates (previously used when frequencies are truncated on some $\cC_{r,R}$) as in our case $r=c_0 2^k$ and $R=C_0 2^j$. We recall that $\varphi$ is the truncation function involved in the Littlewood-Paley decomposition, we denote by $\varphi_1$ another smooth truncation function, with support in a slightly larger annulus than $\varphi$ and equal to $1$ on $supp \varphi$, and by $\cB$ the set:
$$
\cB\overset{def}{=}\{\psi \in \cC_0^\infty (\R_+\times \R^3, \R), \quad \|\psi\|_{L^{\bar{p}}(\R_+, L^{\bar{r}}(\R^3))}\leq 1\},
$$
then following the same classical steps as in \cite{FC3} we get that (we choose for more simplicity to write it only for the third eigenvalue) for any $\beta\geq 1$:
\begin{multline}
 \|\ddk \ddj f\|_{L^p L^r}= \sup_{\psi \in \cB} \int_0^\infty \int_{\R^3} \ddk \ddj f(t,x) \psi(t,x) dx dt\\
 =C \sup_{\psi \in \cB} \int_0^\infty \int_{\R^3} e^{-\nu t|\xi|^2+i\frac{t}{\ee}\frac{|\xi|_F}{F|\xi|}} \widehat{\ddj f_0}(\xi) \varphi_1(2^{-j} \xi) \varphi(2^{-k}|\xi_3|) \hat{\psi}(t,\xi) d\xi dt\\
 \leq C \sup_{\psi \in \cB} \|\ddj f_0\|_{L^2} \left[\int_0^\infty \int_0^\infty \int_{\R^3} K(\nu(t+s), \frac{|t-s|}{\ee},x)\cdot \left( \psi(t)*\bar{\psi}(s)\right)(x) dx ds dt \right]^{\frac12},\\
 \leq C \sup_{\psi \in \cB} \|\ddj f_0\|_{L^2} \left[\int_0^\infty \int_0^\infty \|K(\nu(t+s), \frac{|t-s|}{\ee},.)\|_{L^{\bar{\beta}}} \|\psi(t)*\bar{\psi}(s)\|_{L^\beta} ds dt \right]^{\frac12}
\label{estimTT1}
 \end{multline}
with $K$ defined as follows (we refer to \cite{FC3} for details):
$$
 K(\sigma, \tau,x)=\int_{A_{j,k}} e^{ix\cdot \xi -\sigma |\xi|^2 +i \tau \frac{|\xi|_F}{F|\xi|}} \varphi_1(2^{-j}|\xi|)^2 \varphi(2^{-k}|\xi_3|)^2 d\xi,
$$
where
\begin{equation}
A_{j,k}\overset{def}{=}\{\xi\in \R^3, \quad c_0 2^j \leq |\xi| \leq C_0 2^j \mbox{ and } c_0 2^k \leq |\xi_3| \leq C_0 2^k\}.
\end{equation}
Interpolating the following estimates (we refer to \cite{FC3} Section B.2 for more details), and using as in \cite{FCpochesLp} (section $3.2$) that for all $a,b>0$ and $\theta\in[0,1]$ we have $\min(a,b)\leq a^{1-\theta}b^\theta$:
$$
\begin{cases}
\vspace{0.1cm}
 \|K(\sigma, \tau, .)\|_{L^\infty} \leq C_F e^{-c_0^2\sigma 2^{2j}} 2^{3j} \min(2^{k-j}, \frac{1}{\tau^\frac12 2^{k-j}}),\\
 \|K(\sigma, \tau, .)\|_{L^2} \leq C_F e^{-\frac{c_0^2}{2}\sigma 2^{2j}} 2^{\frac{3j}{2}} 2^{\frac{k-j}{2}},
\end{cases}
$$
we get for any $r\in[2,\infty]$, $\frac{1}{r}=\frac{1-\aa}{\infty}+\frac{\aa}{2}=\frac{\aa}{2}$, and $\theta \in [0,1]$
\begin{multline}
 \|K(\sigma, \tau, .)\|_{L^r} \leq C_F e^{-\frac{c_0^2}{2}\sigma2^{2j}} \left( 2^{3j}\frac{2^{(k-j)(1-2\theta)}}{\tau^{\frac{\theta}{2}}}\right)^{1-\frac{2}{r}} \left(2^{\frac{3j}{2}} 2^{\frac{k-j}{2}}\right)^{\frac{2}{r}}\\
 \leq C_F e^{-\frac{c_0^2}{2}\sigma2^{2j}} 2^{3j(1-\frac{1}{r})}\frac{2^{(k-j)[1-\frac{1}{r}-2\theta(1-\frac{2}{r})]}}{\tau^{\frac{\theta}{2}(1-\frac{2}{r})}}.
 \label{estimK}
\end{multline}
Now we can go back to \eqref{estimTT1}, by the Cauchy-Schwarz inequality, fixing $\beta \geq 1$ so that:
$$
\|\psi(t)*\bar{\psi}(s)\|_{L^\beta} \leq \|\psi(t)\|_{L^{\bar{r}}} \|\psi(s)\|_{L^{\bar{r}}},
$$
that is choosing $\bar{\beta}=\frac{\beta}{\beta-1}=\frac{r}{2}$ (which implies that $r\geq 4$), and using \eqref{estimK}, we obtain that
\begin{multline}
 \|\ddk \ddj f\|_{L^p L^r} \leq C_F \sup_{\psi \in \cB} \|\ddj f_0\|_{L^2} \ee^{\frac{\theta}{4}(1-\frac{4}{r})} 2^{\frac{3j}{2}(1-\frac{2}{r})} 2^{\frac{k-j}{2}(1-\frac{2}{r}-2\theta(1-\frac{4}{r}))}\\
 \times \left[\int_0^\infty \int_0^\infty \frac{h(t)h(s)}{|t-s|^{\frac{\theta}{2}(1-\frac{4}{r})}} ds dt \right]^{\frac12},
\label{estimTT2}
 \end{multline}
with $h(t)=e^{-\frac{c_0^2}{2}\nu t 2^{2j}} \|\psi(t)\|_{L^{\bar{r}}}$. Next we will use the Hardy-Littlewood-Sobolev estimates, that we recall in $\R$ for the convenience of the reader (we refer to \cite{HL, So, Lieb}):
\begin{prop}
 \sl{There exists a constant $C>0$ such that for any function $h_i\in L^{q_i}(\R)$ ($q_i>1$ for $i=1,2$) and any $\aa>0$, with $\frac{1}{q_1}+ \frac{1}{q_2}+\aa =2$, then we have
 $$
 \int_\R \int_\R \frac{h_1(t) h_2(s)}{|t-s|^\aa} dtds \leq C \|h_1\|_{L^{q_1}} \|h_2\|_{L^{q_2}}.
 $$ 
 }
 \label{HL}
\end{prop}
Choosing $h_1=h_2=h \textbf{1}_{\R_+}$, $\aa=\frac{\theta}{2}(1-\frac{4}{r})>0$ and $\frac{1}{q}=1-\frac{\theta}{4}(1-\frac{4}{r})$, we get that
\begin{multline}
 \int_0^\infty \int_0^\infty \frac{h(t)h(s)}{|t-s|^{\frac{\theta}{2}(1-\frac{4}{r})}} ds dt \leq C\|h\|_{L^q}^2 \leq C\left( \|e^{-\frac{c_0^2}{2}\nu 2^{2j}t}\|_{L^m} \|\psi\|_{L^{\bar{p}}L^{\bar{r}}} \right)^2\\
 \leq C\left(\frac{1}{\nu^{\frac{1}{m}}} \left[\frac{2}{mc^2}\right]^{\frac{1}{m}} 2^{-\frac{2j}{m}} \|\psi\|_{L^{\bar{p}}L^{\bar{r}}} \right)^2,
\end{multline}
for $m\in [1,\infty]$ chosen so that $\frac{1}{m}+\frac{1}{\bar{p}}=\frac{1}{q}$, that is:
$$
\frac{1}{m}=\frac{1}{p}-\frac{\theta}{4}(1-\frac{4}{r}).
$$
\begin{rem}
 \sl{Note that this implies the condition $p \leq \frac{4}{\theta} \frac{1}{1-\frac{4}{r}}$.}
\end{rem}
Gathering with \eqref{estimTT2}, we can write that
\begin{multline}
 \|\ddk \ddj f\|_{L^p L^r} \leq C_F \|\ddj f_0\|_{L^2} \ee^{\frac{\theta}{4}(1-\frac{4}{r})} 2^{j(\frac32-\frac{3}{r} -\frac{2}{p} +\frac{\theta}2 (1-\frac{4}{r}))} \frac{2^{\frac{k-j}{2}(1-\frac{2}{r}-2\theta(1-\frac{4}{r}))}}{\nu^{\frac{1}{p}-\frac{\theta}{4}(1-\frac{4}{r})}}\\
 \times \left[\frac{2}{c^2}(\frac{1}{p}-\frac{\theta}{4}(1-\frac{4}{r}))\right]^{\frac{1}{p}-\frac{\theta}{4}(1-\frac{4}{r})}.
\label{estimTT3}
 \end{multline}
It is possible to sum this for $k\leq j+1$ if and only if $1-\frac{2}{r}-2\theta(1-\frac{4}{r})>0$ that is, as $r>4$, when
$$
\theta<\frac12 \frac{1-\frac{2}{r}}{1-\frac{4}{r}}.
$$
Summing over $k$ we obtain that for all $j$,
\begin{equation}
 \|\ddj f\|_{L^p L^r} \leq C_F \frac{C_{p,\theta,r}}{\nu^{\frac{1}{p}-\frac{\theta}{4}(1-\frac{4}{r})}} \ee^{\frac{\theta}{4}(1-\frac{4}{r})} 2^{j(\frac32-\frac{3}{r} -\frac{2}{p} +\frac{\theta}2 (1-\frac{4}{r}))} \|\ddj f_0\|_{L^2},
\end{equation}
which leads to the desired result in the homogeneous case. The inhomogeneous case (i.-e. when $\Fe\neq 0$) easily follows thanks to the Duhamel formula. $\blacksquare$

\subsubsection{Strichartz estimates in the case $\nu\neq\nu'$}

\begin{prop}
 \sl{There exists a constant $C_{F,\omega}>0$ (where $\omega=\frac{\max(\nu, \nu')}{\nu_0}$) such that for any $d\in \R$, $r>4$ and $p<\frac{4}{1-\frac4{r}}$, if $f$ solves \eqref{systdisp} for initial data $f_0$ and external force $\Fe=F^b+F^l$, the three of them with zero divergence and potential vorticity, then for $i=3,4$,
 \begin{multline}
  \||D|^d \mathbb{P}_i \cPrr f\|_{\tilde{L}_t^p\dot{B}_{r, q}^0} \leq C_{F,\omega}\frac{D_{p,r}}{\nu_0^{\frac{1}{p}-\frac14 (1-\frac{4}{r})}} \ee^{\frac14(1-\frac{4}{r})} \frac{\Re^{4-\frac{9}{r}}}{\re^{\frac52+\frac2{p}-\frac6{r}}}\\
  \times \left( \|\cPrr f_0\|_{\dot{B}_{2, q}^d}+ \|\cPrr F^b\|_{L^1 \dot{B}_{2, q}^d} +\frac1{\nu_0^\frac12 \re} \|\cPrr F^l\|_{L^2 \dot{B}_{2, q}^d}\right),
 \end{multline}
 where $D_{p,r}=\max(b_{p,r}, d_{p,r})$ with
$$
\begin{cases}
 b_{p,r}=\left(\frac2{\nu c^2}\right)^{\frac1{p}-\frac14(1-\frac4{r})} \left( \frac1{p}- \frac14(1-\frac4{r}) \right)^{\frac1{p}-\frac14(1-\frac4{r})}\\
 d_{p,r}=2^{\frac1{p}}\left(\frac8{c^2 p}\right)^{\frac{1}{p}-\frac14(1-\frac{4}{r})} \left(\int_0^\infty \frac{e^{-x}}{x^{\frac{p}{4} (1-\frac4{r})}} dx\right)^{\frac1{p}}.
\end{cases}
 $$
\label{Estimdispnu2}
 }
\end{prop}
\begin{rem}
 \sl{We could prove like in the previous section some refined estimate with $\theta\in]0,1]$ (allowing $p\leq \frac4{\theta(1-\frac4{r})}$) but we will only need the case $\theta=1$ and $p$ close to $2$ in this article.}
\end{rem}

\textbf{Proof: } Let us first assume that $\Fe=0$. With the same notations as in the previous section, we get that (see previous section, as well as \cite{FC5, FCpochesLp} for details):
\begin{multline}
 \|\mathbb{P}_i \cPrr f\|_{\tilde{L}^p (\R_+,L^r(\R^3)}=\underset{\psi \in \cB}{\sup} \int_0^\infty \int_{\R^3} \mathbb{P}_i \cPrr f(t,x) \psi(t,x) dx dt\\
 =\underset{\psi \in \cB}{\sup} \int_0^\infty \int_{\R^3} e^{-t\tau (\xi)|\xi|^2+it\frac{|\xi|_F}{\ee F|\xi|}+ \ee t E(\xi, \ee)} \widehat{\mathbb{P}_i \cPrr f_0}(t,\xi) \chi(\frac{|\xi|}{2\Re})\Big(1-\chi(\frac{2|\xi_3|}{\re})\Big)\hat{\psi}(t,\xi) d\xi dt\\
 \leq C \sup_{\psi \in \cB} \|\mathbb{P}_i \cPrr f_0\|_{L^2} \left[\int_0^\infty \int_0^\infty \|L(t,s,\ee,.)\|_{L^\frac{r}2} \|\psi(t)*\bar{\psi}(s)\|_{L^\frac{r}{r-2}} ds dt \right]^{\frac12},
\end{multline}
where
$$
L(t,s,\ee,x)=\int_{\R^3} e^{ix\cdot \xi-(t+s)\tau (\xi)|\xi|^2+i(t-s)\frac{|\xi|_F}{\ee F|\xi|}+ \ee t E(\xi, \ee)+\ee s \bar{E}(\xi, \ee)} \chi(\frac{|\xi|}{2\Re})^2 \Big(1-\chi(\frac{2|\xi_3|}{\re})\Big)^2 d\xi.
$$
Like before, to obtain the $L^{\frac{r}2}$-norm, we will interpolate between $L^2$ and $L^\infty$. It is easy to obtain
$$
\|L(s,t,\ee,.)\|_{L^2} \leq C_F \Re^\frac32 e^{-c^2 \frac{\nu_0}{4}(t+s)\re^2},
$$
and we refer to \cite{FC5, FCpochesLp} where we proved that (there we were working with local in time solutions, and we dropped the exponential):
$$
\|L(s,t,\ee,.)\|_{L^\infty} \leq C_{F,\omega} \frac{\Re^3}{\re^2} \min(1, \frac{\Re^3}{\re^2} \left(\frac{\ee}{|t-s|}\right)^\frac12) e^{-c^2 \frac{\nu_0}{4}(t+s)\re^2}.
$$
so that we obtain for any $\beta\geq 2$:
$$
\|L(s,t,\ee,.)\|_{L^\beta} \leq C_{F,\omega} e^{-c^2 \frac{\nu_0}{4}(t+s)\re^2} \frac{\Re^{6-\frac9{\beta}}}{\re^{4-\frac8{\beta}}} \left(\frac{\ee}{|t-s|}\right)^{\frac12 (1-\frac2{\beta})}.
$$
Thanks to \eqref{estimvp34}, and doing the same as previously, we end-up with ($\beta=\frac{r}2$):
\begin{equation}
 \|\mathbb{P}_i \cPrr f\|_{L^p L^r}=C_{F,\omega} \|\cPrr f_0\|_{L^2} \frac{\Re^{4-\frac9{r}}}{\re^{3-\frac8{r}}} \ee^{\frac14 (1-\frac4{r})} \sup_{\psi \in \cB} \left[\int_0^\infty \int_0^\infty  \frac{g(t) g(s)}{|t-s|^{\frac12(1-\frac{4}{r})}} ds dt \right]^{\frac12},
\end{equation}
with $g(t)=e^{-\frac{c^2}{2}\nu t \re^2} \|\psi(t)\|_{L^{\bar{r}}}$. Using once more Proposition \ref{HL}, we end-up with:
\begin{equation}
  \|\mathbb{P}_i \cPrr f\|_{L^p L^r} \leq C_F\frac{b_{p,r}}{\nu_0^{\frac{1}{p}-\frac14 (1-\frac{4}{r})}} \ee^{\frac14(1-\frac{4}{r})} \frac{\Re^{4-\frac{9}{r}}}{\re^{\frac52+\frac2{p}-\frac6{r}}} \|\cPrr f_0\|_{L^2}.
 \end{equation}
 Then it is easy to deduce the non-homogeneous case with $F^b$ only. Let us now focus on the other external force term, we extend the method from \cite{FC2}. If we denote by $S(t)f_0$ the solution of System \eqref{systdisp} with $\Fe=0$, we have by the Duhamel formula:
\begin{multline}
 \|\int_0^t S(t-t')\cPrr \mathbb{P}_i F^l(t') dt'\|_{L_t^p L^r} =\underset{\psi \in \cB}{\sup} \int_0^\infty \int_{\R^3} \widehat{\mathbb{P}_i \cPrr F^l}(t',\xi)\\
 \times \int_{t'}^\infty e^{-(t-t')\tau (\xi)|\xi|^2+i(t-t')\frac{|\xi|_F}{\ee F|\xi|}+ \ee t E(\xi, \ee)}  \chi(\frac{|\xi|}{2\Re})\Big(1-\chi(\frac{2|\xi_3|}{\re})\Big) \hat{\psi}(t,\xi) dt d\xi dt'\\
  \leq C \sup_{\psi \in \cB} \|\mathbb{P}_i \cPrr F^l\|_{L^2 L^2} \left[\int_0^\infty \int_{t'}^\infty \int_{t'}^\infty \|L(t-t',s-t',\ee,.)\|_{L^\frac{r}2} \|\psi(t)*\bar{\psi}(s)\|_{L^\frac{r}{r-2}} ds dt \right]^{\frac12},\\
  \leq C_{F,\omega} \|\cPrr F^l\|_{L^2 L^2} \frac{\Re^{4-\frac9{r}}}{\re^{3-\frac8{r}}} \ee^{\frac14 (1-\frac4{r})}\\
  \times \sup_{\psi \in \cB} \left[\int_0^\infty \int_0^\infty \int_0^\infty \textbf{1}_{\{t'\leq \min(t,s)\}} \frac{e^{-c^2 \frac{\nu_0}{4}(t+s-2t')\re^2}}{|t-s|^{\frac12(1-\frac4{r})}} \|\psi(t)\|_{L^{\bar{r}}} \|\psi(s)\|_{L^{\bar{r}}} ds dt dt' \right]^{\frac12}.
\end{multline}
Computing the integral in $t'$:
$$
\int_0^{min(s,t)} e^{c^2\frac{\nu_0}{2}t'\re^2} dt'=\frac{2}{\nu_0 \re^2} e^{c^2\frac{\nu_0}{2} \min(t,s)\re^2},
$$
and using the fact that $|t-s|=s+t-2 \min(s,t)$, we get
\begin{multline}
 \|\int_0^t S(t-t')\cPrr \mathbb{P}_i F^l(t') dt'\|_{L_t^p L^r} \leq C_{F,\omega} \|\cPrr F^l\|_{L^2 L^2} \frac{\Re^{4-\frac9{r}}}{\re^{4-\frac8{r}}} \ee^{\frac14 (1-\frac4{r})}\\
  \times \sup_{\psi \in \cB} \left[\int_0^\infty \int_0^\infty \frac{e^{-c^2 \frac{\nu_0}{4}|t-s|\re^2}}{|t-s|^{\frac12(1-\frac4{r})}} \|\psi(t)\|_{L^{\bar{r}}} \|\psi(s)\|_{L^{\bar{r}}} ds dt \right]^{\frac12}.
\end{multline}
Then denoting $k(\tau)=e^{-c^2 \frac{\nu_0}{4}|\tau|\re^2} |\tau|^{-\frac12(1-\frac4{r})}$, we just have to estimate a convolution:
\begin{equation}
 \int_0^\infty \int_0^\infty k(t-s)\|\psi(t)\|_{L^{\bar{r}}} \|\psi(s)\|_{L^{\bar{r}}} ds dt \leq \|k\|_{L^\frac{p}2} \|\psi\|_{L^{\bar{p}} L^{\bar{r}}}^2,
\end{equation}
provided that $p\geq 2$ and $\frac{p}4 (1-\frac4{r})<1$ so that $k\in L^{\frac{p}2}$, whose norm is featured in the constant $d_{p,r}$ and we have
\begin{multline}
  \|\mathbb{P}_i \cPrr f\|_{L^p L^r} \leq C_{F,\omega}\frac{D_{p,r}}{\nu_0^{\frac{1}{p}-\frac14 (1-\frac{4}{r})}} \ee^{\frac14(1-\frac{4}{r})} \frac{\Re^{4-\frac{9}{r}}}{\re^{\frac52+\frac2{p}-\frac6{r}}}\\
  \times \left( \|\cPrr f_0\|_{L^2}+ \|\cPrr F^b\|_{L^1 L^2} +\frac1{\nu_0^\frac12 \re} \|\cPrr F^l\|_{L^2 L^2}\right).
 \end{multline}
Finally, to obtain the announced estimates, we just have to apply this estimates to $\ddj |D|^d f$. $\blacksquare$
\\

\textbf{Aknowledgements :} This work was supported by the ANR project INFAMIE, ANR-15-CE40-0011.

\end{document}